\documentclass[11pt]{article}
\usepackage{amsfonts,amsmath,amsxtra}
\usepackage{latexsym}
\usepackage{amssymb}

\def\hybrid{\topmargin 0pt      \oddsidemargin 0pt
        \headheight 0pt \headsep 0pt
        \textwidth 16.5cm
        \textheight 23cm
        \marginparwidth 0.0in
        \parskip 5pt plus 1pt   \jot = 1.5ex}
\catcode`\@=11
\def\marginnote#1{}
\newcount\hour
\newcount\minute
\newtoks\amorpm
\hour=\time\divide\hour by60 \minute=\time{\multiply\hour by60
\global\advance\minute by-\hour}
\edef\standardtime{{\ifnum\hour<12 \global\amorpm={am}%
        \else\global\amorpm={pm}\advance\hour by-12 \fi
        \ifnum\hour=0 \hour=12 \fi
      \number\hour:\ifnum\minute<10 0\fi\number\minute\the\amorpm}}
\edef\militarytime{\number\hour:\ifnum\minute<10 0\fi\number\minute}

\def\draftlabel#1{{\@bsphack\if@filesw {\let\thepage\relax
   \xdef\@gtempa{\write\@auxout{\string
      \newlabel{#1}{{\@currentlabel}{\thepage}}}}}\@gtempa
   \if@nobreak \ifvmode\nobreak\fi\fi\fi\@esphack}
        \gdef\@eqnlabel{#1}}
\def\@eqnlabel{}
\def\@vacuum{}
\def\draftmarginnote#1{\marginpar{\raggedright\scriptsize\tt#1}}

\def\draft{\oddsidemargin -0.1truein
        \def\@oddfoot{\sl preliminary draft \hfil
        \rm\thepage\hfil\sl\today\quad\militarytime}
        \let\@evenfoot\@oddfoot \overfullrule 3pt
        \let\label=\draftlabel
        \let\marginnote=\draftmarginnote
\def\@eqnnum{{\rm (\theequation)}
\rlap{\kern\marginparsep\tt\@eqnlabel}%
\global\let\@eqnlabel\@vacuum}  }
\newfont{\Bbbb}{msbm7 scaled 1\@ptsize00}
\newcommand{\zs}{\raise-1pt\hbox{$\mbox{\Bbbb Z}$}}

scaled 1\@ptsize00
\font\sevenmsa=msam6 %scaled 1\@ptsize00
\newfam\msafam
\textfont\msafam=\sevenmsa
\def\hexnumber@#1{\ifnum#1<10 \number#1\else
\ifnum#1=10 A\else\ifnum#1=11 B\else\ifnum#1=12 C\else \ifnum#1=13
D\else\ifnum#1=14 E\else\ifnum#1=15 F\fi\fi\fi\fi\fi\fi\fi}
\def\msa@{\hexnumber@\msafam}
\def\llcorner{\delimiter"4\msa@78\msa@78 }
\def\lrcorner{\delimiter"5\msa@79\msa@79 }
\mathchardef\blacktriangleright="3\msa@49
\mathchardef\blacktriangleleft="3\msa@4A \font\tenmsb=msbm10 scaled
1\@ptsize00
\newfam\msbfam
\textfont\msbfam=\tenmsb \scriptfont\msbfam=\tenmsb
\newdimen\Squaresize \Squaresize=14pt
\newdimen\Thickness \Thickness=0.5pt

\def\Square#1{\hbox{\vrule width \Thickness
   \vbox to \Squaresize{\hrule height \Thickness\vss
      \hbox to \Squaresize{\hss#1\hss}
   \vss\hrule height\Thickness}
\unskip\vrule width \Thickness} \kern-\Thickness}

\def\Vsquare#1{\vbox{\Square{$#1$}}\kern-\Thickness}

\def\numberbysection{\@addtoreset{equation}{section}
        \def\theequation{\thesection.\arabic{equation}}}
\numberbysection

\renewcommand{\theequation}{\thesection.\arabic{equation}}
\def\titlepage{\@restonecolfalse\if@twocolumn\@restonecoltrue\onecolumn
     \else \newpage \fi \thispagestyle{empty}\c@page\z@
        \def\thefootnote{\fnsymbol{footnote}} }

\def\endtitlepage{\if@restonecol\twocolumn \else  \fi
        \def\thefootnote{\arabic{footnote}}
        \setcounter{footnote}{0}}  %\c@footnote\z@ }
\relax

\hybrid
\makeatletter
\newdimen\normalarrayskip            % skip between lines
\newdimen\minarrayskip               % minimal skip between lines
\normalarrayskip\baselineskip \minarrayskip\jot
\newif\ifold             \oldtrue            \def\new{\oldfalse}
\def\arraymode{\ifold\relax\else\displaystyle\fi}%mode of array enrties
\def\eqnumphantom{\phantom{(\theequation)}} % ight phantom in eqnarray
\def\@arrayskip{\ifold\baselineskip\z@\lineskip\z@
     \else
     \baselineskip\minarrayskip\lineskip1\baselineskip\fi}

\def\@arrayclassz{\ifcase \@lastchclass \@acolampacol \or
\@ampacol \or \or \or \@addamp \or
   \@acolampacol \or \@firstampfalse \@acol \fi
\edef\@preamble{\@preamble
  \ifcase \@chnum
     \hfil$\relax\arraymode\@sharp$\hfil
     \or $\relax\arraymode\@sharp$\hfil
     \or \hfil$\relax\arraymode\@sharp$\fi}}

\def\@array[#1]#2{\setbox\@arstrutbox=\hbox{\vrule
     height\arraystretch \ht\strutbox
     depth\arraystretch \dp\strutbox
width\z@}\@mkpream{#2}\edef\@preamble{\halign \noexpand\@halignto
\bgroup \tabskip\z@ \@arstrut \@preamble \tabskip\z@ \cr}%
\let\@startpbox\@@startpbox \let\@endpbox\@@endpbox
    \if #1t\vtop \else \if#1b\vbox \else \vcenter \fi\fi
  \bgroup \let\par\relax
  \let\@sharp##\let\protect\relax
  \@arrayskip\@preamble}
\def\eqnarray{\stepcounter{equation}%
              \let\@currentlabel=\theequation
              \global\@eqnswtrue
              \global\@eqcnt\z@
              \tabskip\@centering              %formulae  centering
              \let\\=\@eqncr
              $$%
            \halign to \displaywidth  \bgroup
             \eqnumphantom \@eqnsel
      \hskip\@centering                               %right tab%
    $\displaystyle  \tabskip\z@ {##}$%
    &\global\@eqcnt\@ne \hskip 2\arraycolsep
         $ \displaystyle  \arraymode{##}$\hfil
    &\global\@eqcnt\tw@ \hskip 2\arraycolsep
         $\displaystyle\tabskip\z@{##}$\hfil
         \tabskip\@centering
    &{##}\tabskip\z@\cr}
\makeatother
\newcommand{\CC}{{\mathbb{C}}}
\def\IC{\mathbb{C}}

\def\IQ{\mathbb{Q}}
\def\IR{\mathbb{R}}
\def\IZ{\mathbb{Z}}
\def\CC {\mathcal{C}}
\def\CD {\mathcal{D}}

\def\CH {\mathcal{H}}

\def\CL {\mathcal{L}}
\def\CM {\mathcal{M}}

\def\CQ {\mathcal{Q}}

\def\a {{\alpha}}

\def\g {{\gamma}}
\def\s {{\sigma}}

\def\la{\lambda}

\def\e{\epsilon}
\def\pr {\partial}

\def\nn{\nonumber}

\def\frak{\mathfrak}
\newtheorem{te}{Theorem}[section]%Usage:\begin{te}Statement\end{te}
\newtheorem{de}{Definition}[section]
\newtheorem{prop}{Proposition}[section]
\newtheorem{lem}{Lemma}[section]
\newtheorem{ex}{Example}[section]
\newtheorem{rem}{Remark}[section]

\newcommand{\proof}{\noindent {\it Proof}.\,\,}
\newcommand\bqa{\begin{eqnarray}}
\newcommand\eqa{\end{eqnarray}}
\def\be{\begin{eqnarray}\new\begin{array}{cc}}
\def\ee{\end{array}\end{eqnarray}}
\def\nn{\nonumber}

\def\beq{\begin{equation}}
\def\eeq{\end{equation}}
\def\bse{\begin{subequations}}                %%%SUBEQUATIONS
\def\ese{\end{subequations}}
\def\bp{\begin{pmatrix}}
\def\ep{\end{pmatrix}}

\def\h{\hbar}
\def\i{\imath}

\def\stack#1#2{\raise0.7pt\hbox{$\mathrel{\mathop{#2}\limits^{#1}}$}}
\def\tr{\triangleright}
\def\tl{\triangleleft}
\def\sem{\mathsurround=0pt \raise1pt
\hbox{$\scriptscriptstyle>\!\!$}\:\!\!\tl}
\def\mes{\mathsurround=0pt \tr\!\:\!\raise0.8pt
\hbox{$\scriptscriptstyle\!\!<$}\,}
\def\]{\mathsurround=0pt ]\raise-2pt\hbox{$_\ast$}}
\def\<{\langle}
\def\>{\rangle}

\def\CQ{{\cal Q}}
\def\frak{\mathfrak}

\def\CH{\mathcal{H}}

\def\we{\raise-1pt\hbox{$\,\stackrel{\wedge}{,}\,$}}
\def\tr{{\rm tr}\,}

\def\pr {\partial}

\newcounter{pac}[section]

\newcounter{pacc}[subsection]

  0

\begin{document}

\title{\bf Baxter operator formalism for Macdonald polynomials}
\author{Anton Gerasimov, Dimitri Lebedev and Sergey Oblezin}
\date{}
\maketitle

\renewcommand{\abstractname}{}

\begin{abstract}
\noindent {\bf Abstract}. We develop basic constructions of the
Baxter operator formalism for the Macdonald polynomials associated
with root systems of type A. Precisely we construct a dual pair of
mutually commuting Baxter operators such that the Macdonald
polynomials are their common eigenfunctions. The dual pair of Baxter
operators is closely related to the dual pair of recursive operators
for Macdonald polynomials leading  to various families of their
integral representations. We also construct the Baxter operator
formalism for the $q$-deformed ${\frak gl}_{\ell+1}$-Whittaker
functions and the Jack polynomials obtained by degenerations of the
Macdonald polynomials associated with the type $A_{\ell}$ root
system. This note provides a generalization of our previous results
on the Baxter operator formalism  for the Whittaker functions.  It
was demonstrated previously that  Baxter operator formalism for the
Whittaker functions has deep connections with representation theory.
In particular the Baxter operators  should be considered as elements
of appropriate spherical Hecke algebras and their eigenvalues are
identified with local Archimedean $L$-factors associated with
admissible representations  of reductive groups over $\IR$. We
expect that the Baxter operator formalism for the Macdonald polynomials
has an interpretation in  representation theory of
 higher-dimensional arithmetic fields.

\end{abstract}
\vspace{5 mm}

\section*{Introduction}

A new class of operators acting on eigenfunctions of quantum
integrable systems was introduced by Baxter to provide a solution of
a class of integrable models \cite{Ba}. These operators commute with
quantum Hamiltonians of a quantum integrable system and satisfy
difference/differential  equations with coefficients expressed
through quantum Hamiltonians.  The Baxter operators were constructed
for many integrable models including periodic Toda chains \cite{PG}.

 In \cite{GLO1} we introduce
Baxter operators for non-periodic $\mathfrak{gl}_{\ell+1}$-Toda
chains given by  one-parameter families  of integral operators. We
also define  dual Baxter operators   acting on the spectral
variables of Toda chain eigenfunctions. The dual pair of Baxter
operators  enters a canonical construction of a pair of  recursive
operators relating eigenfunctions of the
$\mathfrak{gl}_{\ell+1}$- and $\mathfrak{gl}_{\ell}$- Toda chains.
This gives rise to various families of integral representations for
eigenfunctions. Hence the Baxter operator formalism consisting of a
pair of dual Baxter operators, a pair of dual recursive operators
provides a complete solution of the Toda chains.

One can expect that the Baxter operator formalism can be constructed
for a wide class of quantum integrable systems. Note that
$\mathfrak{gl}_{\ell+1}$-Toda chains can be considered as a
degeneration of the quantum integrable  system constructed by
Ruijsenaars \cite{Ru} and Macdonald \cite{M}. The corresponding
quantum Hamiltonians are given by mutually commuting difference
operators and their common polynomial eigenfunctions are given by
the Macdonald polynomials. In this note we construct Baxter operator
formalism for the  Macdonald-Ruijsenaars integrable
 system. This includes a dual pair of Baxter operators,
 a dual pair of recursive operators and
various families of explicit iterative expressions for the Macdonald
polynomials, known and new ones.  We also  describe  Baxter operator
formalism for specializations of the  Macdonald polynomials given by
class one $q$-deformed $\mathfrak{gl}_{\ell+1}$-Whittaker functions
\cite{GLO2} and Jack's polynomials. Due to the results of
\cite{GLO8} the Baxter operator formalism for the standard Whittaker
functions \cite{GLO1} can be recovered from the Baxter operator
formalism for $q$-deformed $\mathfrak{gl}_{\ell+1}$-Whittaker
polynomials in the limit $q\to1$.

One should stress that  the Baxter operators associated with the
Whittaker functions have a surprising relation  with number theory
and representation theory \cite{GLO1}. Recall that the
eigenfunctions of $\mathfrak{gl}_{\ell+1}$-Toda chains can be
identified with
 particular matrix elements in the principal series
representations of $GL_{\ell+1}(\IR)$ and thus providing
generalizations of the classical  Whittaker functions corresponding
to $SL_2(\IR)$. In \cite{GLO1} we argue that for
$\mathfrak{gl}_{\ell+1}$-Toda chain the Baxter operator should be
considered as a generating function of elements of spherical Hecke
algebra associated  with the maximal compact  subgroup of
$GL_{\ell+1}(\IR)$. Furthermore the
$\mathfrak{gl}_{\ell+1}$-Whittaker functions are eigenfunctions  of
the Baxter operators with the eigenvalues given by local Archimedean
$L$-factors of the corresponding principle series representations of
$GL_{\ell+1}(\IR)$.

This interpretation leads to  establishing  a deep relation between
topological field theories and Archimedean algebraic  geometry
\cite{GLO5}, \cite{GLO6}, \cite{GLO7}. The construction of the
Baxter operator formalism for the Macdonald-Ruijsenaars integrable
system allows to define a $(q,t)$-generalization  of the local
Archimedean $L$-factors associated with principle series
representations of $GL_{\ell+1}$ ($q$-generalization of local
$L$-factors was introduced previously in \cite{GLO2}). One can
expect that  these generalized $L$-factors shall be related with
principal series representations of  loop groups associated with
$GL_{\ell+1}$. Taking into account the results of \cite{GLO5},
\cite{GLO6}, \cite{GLO7}  one should look for a higher dimensional
topological field theory interpretation of the Macdonald polynomials
and the associated Baxter operator formalism. We are going to
discuss this interpretation in the future publications.

Finally we would like to point out  that many of the constructions
 of this note are simple reformulations of the results of Macdonald
 \cite{M}. However  we feel that  establishing the  direct relation
of the results of \cite{M} with Baxter operator formalism  might be useful.

{\em Acknowledgments}:  This work is partially supported by the RFBR
grant 12-01-00894-a. AG was partly supported by Science Foundation
Ireland grant. The research of AG and DL was also partially
supported by QGM (Centre for Quantum Geometry of Moduli Spaces)
funded by the Danish National Research Foundation. The research of
SO was partially supported by P. Deligne's 2004 Balzan Prize in
Mathematics.

\section{Preliminaries on symmetric polynomials}

In this Section  we collect basic facts on the Macdonald symmetric
polynomials and their degenerate versions given by Jack polynomials
and class one $q$-Whittaker functions. For details on Macdonald and
Jack's polynomials see \cite{M}; for  class one $q$-deformed
Whittaker functions see \cite{GLO2}, \cite{GLO3}, \cite{GLO4}.

\subsection{Macdonald symmetric polynomials}

Let $\IQ(q,t)$ be a field of rational functions in
variables $q$ and $t$. Define the following  $(q,t)$-analog of the
classical $\Gamma$-function (see Appendix for its basic properties):
 \be
  \Gamma_{q,t}(x)\,=\,\frac{(tx;\,q)_{\infty}}{(x;\,q)_{\infty}}\,,\qquad
 (x;\,q)_{\infty}\,=\,\prod_{j=0}^{\infty}\bigl(1-xq^j\bigr).
\ee Let $\Lambda_{q,t}$ be the  graded $\IQ(q,t)$-algebra of
symmetric polynomials  of variables $x_1,x_2,\ldots $ of degree one
 \be
  \Lambda_{q,t}\,=\,\bigoplus_{n\geq0}\,\Lambda_{q,t}^{(n)}\,,
 \ee
where  $\Lambda_{q,t}^{(n)}$ is the homogeneous component of
$\Lambda_{q,t}$ of degree $n$. There are various convenient bases in
the space of symmetric polynomials in variables $x_1,\ldots
,x_{\ell+1}$ enumerated by  partitions
$\la=(\la_1\geq\ldots\geq\la_{\ell+1})$, $\la_i\in \IZ_+$.
Particularly, the elements of the bases of monomial symmetric
functions $m_{\la}(x)$ are given by  sums of all distinct monomials
obtained from $x^{\la}=x_1^{\la_1}\ldots x_{\ell+1}^{\la_{\ell+1}}$
by permutations of $x_1,\ldots,x_{\ell+1}$. Let us denote
$p_{n}(x)\,:=\,m_{(n)}$ the symmetric polynomial for the partition
$(n)=(n,0,\ldots,0)$. The bases of power series symmetric
  polynomials consists of the polynomials
$p_{\la}(x)=p_{\la_1}(x)\cdot\ldots\cdot p_{\la_{\ell+1}}(x)$. Equip the
space $\Lambda_{q,t}^{(\ell+1)}$ with a scalar product $\<\,,\,\>$
defined by
 \be\label{SP0}
 \hspace{-5mm}
  \<p_{\la},\,p_{\mu}\>_{q,\,t}\,
  =\,\delta_{\la\mu}\,z_{\la}
  \prod_{i=1}^{\ell+1}\frac{1-q^{\la_i}}{1-t^{\la_i}}\,,\qquad
  z_{\la}\,=\,\prod_{n\geq1}n^{m_n}m_n!\,,\qquad
  m_n\,=\,\bigl|\{k:\,\la_k=n\}\bigr|\,.
 \ee
Macdonald introduced a bases $\{P_{\la}(x)=P_{\la}(x;\,q,t)\}$ of
symmetric polynomials over $\IQ(q,t)$ enumerated by partitions $\la$
such that
 \be\label{rel0}
  P_{\la}(x)\,=\,\sum_{\mu\leq\la}u_{\la\mu}\,m_{\mu}\,,
 \hspace{1.5cm}
  u_{\la\la}=1\,,
 \ee
and
 \be
  \<P_{\la},\,P_{\mu}\>_{q,\,t}\,=\,0,\,
 \hspace{1.5cm}
  \la\neq\mu\,.
 \ee
In the above formula $\leq$ denotes the natural ordering:
$$
 \la\leq\mu\quad\Longleftrightarrow\quad
 \la_1+\ldots+\la_i\,\leq\,\mu_1+\ldots+\mu_i\,,\quad i\geq0\,.
$$

The  relation \eqref{rel0} is invertible and thus the Macdonald
polynomials $P_{\la}(x)$ provide a bases in
$\Lambda_{q,t}^{(\ell+1)}$. The inverse norms of the Macdonald
polynomials are  given by
 \be\label{Norm}
  b_{\la}\,:=\<P_{\la},\,P_{\la}\>_{q,\,t}^{-1}\,
  =\,\prod_{(i,j)\in\la}
  \frac{1-t^{\la^\top_j+1-i}q^{\la_i-j}}{1-t^{\la^\top_j-i}q^{\la_i+1-j}}\,,
 \ee
where the product is over the boxes $(i,j)$ in the Young diagram
 attached to partition $\la$, and $\la^\top$ denotes the
conjugate partition. In particular, one has
 \be
  b_{(n)}\,
  =\,\prod_{i=1}^n\frac{1-tq^{n-i}}{1-q^{n+1-i}}\,
  =\,\frac{\Gamma_{q,tq^{-1}}(q)}{\Gamma_{q,tq^{-1}}(q^{n+1})}\,.
 \ee

The Macdonald polynomials $P_{\la}(x)$ can be also characterized as
common  eigenfunctions of the following  set of mutually commuting
difference operators \cite{M}, \cite{Ru}  acting in
$\Lambda_{q,t}^{(\ell+1)}$:
 \be\label{MacdonaldHamiltonians}
  M_r\,
  =\,t^{r(r-1)/2}\sum_{I_r}\,
  \prod_{i\in I_r\atop j\notin I_r}\,
  \frac{tx_i-x_j}{x_i-x_j}\,T_{I_r}\,,
 \hspace{1.5cm}
  T_{I_r}\,=\,\prod_{i\in I_r}T_{q,\,x_i}\,,
 \ee
where the sum goes over all $r$-element subsets $I_r$ of
$(1,2,\ldots,\ell+1)$ and
$$
 T_{q,\,x_i}\cdot f(x_1,\ldots,x_{\ell+1})\,\,
 =\,\,f(x_1,\ldots,qx_i,\ldots,x_{\ell+1}).
$$
The operators  $M_r,\,r=1,\ldots,\ell+1$ are self-adjoint with respect
to the scalar product \eqref{SP0} and
the eigenvalues of  $M_r$ acting on $P_{\la}(x)$ are
given by  the elementary symmetric functions
 \be\label{MacdonaldEigenvalues1}
  \chi_r(y_1,\ldots, y_{\ell+1})\,
  =\,\sum_{1\leq i_1<\ldots<i_r\leq \ell+1}y_{i_1}\cdot\ldots\cdot
  y_{i_r}\,,\qquad y_i=t^{\varrho_i}q^{\la_i},
 \ee
where $\varrho_i=\ell+1-i$.
The eigenfunction property of the Macdonald polynomials can be
 succinctly described by the relation
 \be\label{MacdonaldEigenvalues}
  \CM_{\ell+1}(X)\cdot P_{\la}(x)\,
  =\,c_{\ell+1}(\la;\,X)\,P_{\la}(x)\,,\\
  c_{\ell+1}(\la;\,X)\,
  =\,(1-t)^{-(\ell+1)}\prod_{i=1}^{\ell+1}\bigl(1+t^{\varrho_i}q^{\la_i}X\,\bigr)\,,
 \ee
where $\CM_{\ell+1}(X)$ is a generating function of the difference
operators \eqref{MacdonaldHamiltonians}:
 \be\label{MacdonaldGenFunction}
  \CM_{\ell+1}(X)\,\,
  =\,\,(1-t)^{-(\ell+1)}\Big(1+\sum_{r=1}^{\ell+1}\,X^r\,M_r\Big)\,.
 \ee

Assume now that $q\in\IC$ and $|q|<1$,
so that infinite product $(z;\,q)_{\infty}$ converges for all
$z\in\IC$. Following Macdonald define a new  scalar product
 \be\label{AnotherScalarProd}
  \<a,\,b\>'_{q,\,t}\,=\,\frac{1}{(\ell+1)!}\int_Td^{\times}z\,\,
  a(z)\,b(z^{-1})\,\,\Delta(z;\,q,t)\,,\\
  \Delta(z;\,q,t)\,
  =\,\prod_{i,j=1\atop i\neq j}^{\ell+1}
  \frac{1}{\Gamma_{q,t}(z_iz_j^{-1})}\,,
 \ee
where $a(z)$ and $b(z)$ are Laurent polynomials  and
$$
 T\,
 =\,\bigl\{z=(z_1,\ldots,z_{\ell+1})\in\IC^{\ell+1}:\,
 |z_i|=1,\,i=1,\ldots,\ell+1\bigr\}
$$
is the $(\ell+1)$-dimensional  torus, with the Haar measure
$d^{\times}z=\prod_{i=1}^{\ell+1}
(2\pi\i)^{-1}d\log z_i$. The polynomials $P_{\la}$ are pairwise
orthogonal with respect to the new scalar product $\<,\>'$ with the
norms given by
 \be\label{AnotherNorm}
  \<P_{\la},\,P_{\la}\>'_{q,\,t}\,
  =\,\prod_{i,j=1\atop i<j}^{\ell+1}
  \frac{\Gamma_{q,tq^{-1}}(t^{j-i-1}q^{\la_i-\la_j+1})}
  {\Gamma_{q,tq^{-1}}(t^{j-i}q^{\la_i-\la_j+1})}\,.
 \ee

Let us recall the properties of the Macdonald polynomials that will
play essential role in the following (see \cite{M} for the proofs).

\begin{te} Consider two sets $x=(x_1,\ldots,x_n)$ and
$y=(y_1,\ldots,y_m)$ of variables. Let
 \be\label{CauchyLittlewood}
  \Pi_{n,\,m}(x,y):=\sum_{\la\in Y_{n,m}}\,b_{\la}\,P_{\la}(x)\,P_{\la}(y),
 \ee
where summation goes over a set $Y_{m,n}$ of the partitions of length
$\min(m,n)$ and $b_{\lambda}$ are given by \eqref{Norm}. Then the following identity
holds
\be\label{CauchyKer}
\Pi_{n,\,m}(x,y)\,
=\,\prod_{i=1}^n\prod_{j=1}^m
\Gamma_{q,t}(x_iy_j).
\ee
\end{te}
We will say that partitions $\mu,\lambda$ are interlaced if
$\mu_1\geq\lambda_1\geq\ldots\geq\mu_{\ell+1}\geq\lambda_{\ell+1}$.
In the sequel we shall use the following abbreviation for interlaced
partitions: $\mu_i\geq\lambda_i\geq\mu_{i+1}$.
\begin{te} Let $P_{(n)}(x)$ be the Macdonald polynomial
  corresponding to the partition $(n)=(n,0,\ldots,0)$. Then
the following product decomposition holds:
 \be\label{PieriRules1}
  P_{(n)}(x)\times P_{\la}(x)\,\,
  =\,\,\,b_{(n)}^{-1}\sum_{\mu_i\geq\la_i\geq\mu_{i+1}\atop |\mu|-|\la|=n}
  \varphi_{\mu/\la}\,\,P_{\mu}(x)\,,
 \ee
where
 \be\label{PieriRules2}
  \varphi_{\mu/\la}\,
  =\,\prod_{i,j=1\atop i\leq j}^{\ell+1}
  \frac{\Gamma_{q,tq^{-1}}\bigl(t^{j-i}q^{\mu_i-\mu_j+1}\bigr)}
  {\Gamma_{q,tq^{-1}}\bigl(t^{j-i}q^{\mu_i-\la_j+1}\bigr)}\,
  \frac{\Gamma_{q,tq^{-1}}\bigl(t^{j-i}q^{\la_i-\la_{j+1}+1}\bigr)}
  {\Gamma_{q,tq^{-1}}\bigl(t^{j-i}q^{\la_i-\mu_{j+1}+1}\bigr)}\,,
 \ee
and we omit in the product \eqref{PieriRules2} the factors
depending on $\la_{\ell+2}$ and $\mu_{\ell+2}$.
\end{te}

The Macdonald polynomials possess a remarkable self-duality property
discovered by Koornwinder (see \cite{M} and references
therein). Let us introduce modified Macdonald polynomials:
 \be
  \Phi_{\la}(x;\,q,t)\,:=\,\,t^{\rho(\la)}
  \prod_{a,b=1\atop a<b}^{\ell+1}
  \Gamma_{q,t}(t^{b-a}q^{\la_a-\la_b})\times
  P_{\la}(x;\,q,t)\,,
 \ee
where $\rho(x)=\sum\limits_{i=1}^{\ell+1}\rho_ix_i$, $\rho_i=\varrho_i-\ell/2$.
Then  for any partitions $\lambda$ and $\mu$ the following duality relation
holds:
 \be\label{SelfDuality}
  \Phi_{\la}\bigl(q^{\mu-k\rho};\,q,\,q^{-k}\bigr)\,
  =\,\Phi_{\mu}\bigl(q^{\la-k\rho};\,q,\,q^{-k}\bigr)\,.
 \ee
This duality naturally leads to the set of mutually commuting
operators acting in the space of functions on the set of partitions
$\lambda=(\lambda_1,\ldots ,\lambda_{\ell+1})\in \IZ_+^{\ell+1}$.
The following Theorem was proved in \cite{GLO4}.

\begin{te}[\cite{GLO4}]
A set of mutually commuting difference operators
 \be\label{MacdonaldDualHamiltonians2}
  M_r^{\vee}(\la)\,
  =\,t^{r\ell/2}\sum_{I_r}\prod_{i\in I_r\atop{j\notin I_r\atop j<i}}
  \frac{1-t^{i-j+1}q^{\la_j-\la_i-1}}{1-t^{i-j}q^{\la_j-\la_i-1}}
  \frac{1-t^{i-j-1}q^{\la_j-\la_i}}{1-t^{i-j}q^{\la_j-\la_i}}\,
  T^{\vee}_{I_r}\,,\\
  T_{q,\,q^{\la_i}}\cdot f_{(\la_1,\ldots,\la_{\ell+1})}\,
  =\,f_{(\la_1,\ldots,\la_i+1,\,\ldots,\la_{\ell+1})}\,,
 \hspace{1.5cm}
  T^{\vee}_{I_r}=\prod_{i\in I_r}T_{q,\,q^{\la_i}}\,,
 \ee
acts in the space of functions $f_{\la}$ labeled by partitions
$\lambda=(\lambda_1\geq\ldots\geq\lambda_{\ell+1})\in
\IZ_+^{\ell+1}$. The Macdonald polynomials $P_\lambda(x)$  as
functions of the variables $\lambda$ are common eigenfunctions of
the difference operators \eqref{MacdonaldDualHamiltonians2}:
 \be\label{DualEigenvalues}
  \CM_{\ell+1}^{\vee}(X)\cdot P_{\la}(x;\,q,t)\,\,
  =\,\,c^{\vee}_{\ell+1}(t^{\ell/2}x;\,X)\,P_{\la}(x;\,q,t)\,,\\
  c^{\vee}_{\ell+1}(x;\,X)\,
  =\,\prod_{i=1}^{\ell+1}\bigl(1+Xx_i\bigr)\,,
\ee
where
 \be\label{MacdonaldDualGenFunction}
  \CM_{\ell+1}^{\vee}(X)\,
  =\,\sum_{r=0}^{\ell+1}\,X^r\,M^{\vee}_r\,,
 \hspace{1.5cm}
  M^{\vee}_0:=1\,,
\ee
\end{te}

\subsection{Class one $q$-deformed Whittaker functions}

Let $\Lambda_q$ be the algebra of symmetric functions in variables
$x_1,x_2,\ldots$ over the field of rational functions in variable
$q$ convergent in the domain $|q|<1$. Let $\Lambda_q^{(\ell+1)}$ be
the homogeneous component $\Lambda_q$ of degree $\ell+1$. Consider a
pair of the scalar products on $\Lambda_q^{(\ell+1)}$. The first one
is defined in terms of power series symmetric polynomials $p_{\la}$
as follows:
 \be\label{Wsp1}
  \<p_{\la},\,p_{\mu}\>_q\,
  =\,\delta_{\la\mu}\,z_{\la}\,
  \prod_{i=1}^{\ell+1}(1-q^{\la_i})\,.
 \ee
The second  scalar product on the space  $\Lambda_q$ is defined by
 \be\label{Wsp2}
  \<a(z),\,b(z)\>'_q\,
  =\,\frac{1}{(\ell+1)!}\int_T\!d^{\times}\!z\,\,
  a(z)\,b(z^{-1})\,\Delta^{\vee}_q(z)\,,\\
  \Delta^{\vee}_q(z)\,=\,\prod_{i,j=1\atop i\neq j}^{\ell+1}
  \frac{1}{\Gamma_q(z_i^{-1}z_j)}\,,\qquad
\Gamma_q(x)=\prod_{j=0}^{\infty}\frac{1}{1-q^jx},
\ee
and the notations of  \eqref{AnotherScalarProd} are used.

\begin{lem}
 The polynomials $P_{\lambda}(x;q):=P_{\lambda}(x;q,t=0)$
satisfy the relations
 \be\label{rel1}
  P_{\la}(x)\,=\,\sum_{\mu\leq\la}u_{\la\mu}\,m_{\mu}\,,
 \hspace{1.5cm}
  u_{\la\la}=1\,,
 \ee
 \be\label{orth1}
  \<P_{\la},\,P_{\mu}\>_{q,\,t}\,=\,0,\,
 \hspace{1.5cm}
  \la\neq\mu\,.
 \ee
and thus  define  a bases in $\Lambda_q$.
\end{lem}

\proof Directly follows by specialization $t=0$
from the properties of the  Macdonald polynomials. $\Box$

Let us define the normalized symmetric polynomials in variables
$x_1,\ldots,x_{\ell+1}$, labeled by partitions
$\lambda=(\lambda_1\geq\ldots\geq\lambda_{\ell+1})\in\IZ_+^{\ell+1}$
as follows
 \be\label{WhitPol}
  P^{qW}_{\lambda}(x)\,
  =\Delta_q^{-1}(\lambda)P_{\lambda}(x;q,t=0)\,,
 \hspace{1.2cm}
  \Delta_q(\la)\,=\,\prod_{i=1}^{\ell}(\la_i-\la_{i+1})_q!\,,
 \ee
where $(n)_q!=\prod_{i=1}^n(1-q^{n+1-i})$. In the following we will
call $P^{qW}_{\lambda}(x)$  the $q$-Whittaker polynomials.
These polynomials were introduced in \cite{GLO2} as class one
$q$-deformed $\mathfrak{gl}_{\ell+1}$-Whittaker functions
$\Psi^{(q)}_x(\lambda)=P^{qW}_{\lambda}(x)$.

The $q$-Whittaker polynomials $P^{qW}_{\lambda}(x)$ are
  orthogonal with respect to both scalar products \eqref{Wsp1}, \eqref{Wsp2}  and
are normalized as follows:
 \be\label{qWhittakerNorm}
  \<P^{qW}_{\la},\,P^{qW}_{\la}\>_q\,
  =\,\frac{(\la_{\ell+1})_q!}{\Delta_q(\la)}\,,
 \hspace{1.5cm}
  \<P^{qW}_{\la},\,P^{qW}_{\la}\>'_q\,
  =\,\frac{1}{\Delta_q(\la)}\Big(\frac{1}{\Gamma_q(q)}\Big)^{\ell}\,.
 \ee

\begin{te}[\cite{GLO4}]  Let $H_1,\ldots ,H_{\ell+1}$, and
$H^{\vee}_1,\ldots ,H^{\vee}_{\ell+1}$ be  difference operators
acting in the space of functions on $\IR^{\ell+1}\times
\IZ_+^{\ell+1}$
 \be\label{qWhittakerDualHamiltonians}
  H^{\vee}_r\,
  =\,\sum_{I_r}\prod_{i\in I_r\atop j\notin I_r}
  \frac{x_j}{x_j-x_i}\,T_{I_r}\,,
 \hspace{1.5cm}
  T_{I_r}\,=\,\prod_{i\in I_r}T_{q,\,x_i}\,,\qquad r=1,\ldots,\ell+1,
 \ee
 \be\label{qWhittakerHamiltonians}
  H_r\,=\,\sum_{I_r}\prod_{k=1}^r
  \bigl(1-q^{\la_{i_k}-\la_{i_k+1}+1}\bigr)^{1-\delta_{i_{k+1}-i_k,\,1}}\,
  T^{\vee}_{I_r}\,,
 \hspace{1.5cm}
  T^{\vee}_{I_r}\,=\,\prod_{i\in I_r}T_{q,\,q^{\la_i}}\,,
 \ee
for $r=1,\ldots,\ell+1$, where $i_{r+1}:=\ell+2$ is  assumed.

These  operators are mutually commutative
and the  $q$-Whittaker polynomials solve the following dual
pair of eigenfunction problems:
 \be
  H_r^{\vee}\cdot P^{qW}_\lambda(x)
  =\,q^{\la_{\ell+2-r}+\ldots+\la_{\ell+1}}\,P^{qW}_\la(x),\qquad
  r=1,\ldots,\ell+1\,.
 \ee
and
 \be
  H_r\cdot P^{qW}_{\la}(x)\,
  =\,\chi_r(x)\,P^{qW}_{\la}(x),\,
 \hspace{1.5cm}
  \chi_r(x)\,=\,\sum_{I_r}x_{i_1}\cdots x_{i_r}\,,
 \ee
for $r=1,\ldots,\ell+1$.
\end{te}
For the generating function $D_{\ell+1}$ of the operators $H_r$
 \be
  D_{\ell+1}(X)\,=\,\sum_{r=0}^{\ell+1}X^rH_r\,,
 \hspace{1.5cm}
  H_0:=1\,,
 \ee
the following relation holds
 \be\label{qWhitEigenvalues}
  D_{\ell+1}(X)\cdot P^{qW}_{\la}(x)\,
  =\,{\bf c}^q_{\ell+1}(x;\,X)\, P^{qW}_{\la}(x)\,,\\
  {\bf c}^q_{\ell+1}(x;\,X)\,
  =\,\prod_{i=1}^{\ell+1}(1+X\,x_i)\,.
 \ee

The set of operators \eqref{qWhittakerHamiltonians}  define
$q$-deformed Toda chain  Hamiltonians and
\eqref{qWhittakerDualHamiltonians} provide a set of mutually
commuting difference dual Toda chain Hamiltonians introduced  in
\cite{GLO4}.

Using the relation between $q$-Whittaker polynomials $P^{qW}_\la(x)$
and Macdonald polynomials $P_\la(x)$ one can infer an analog of the
Pieri formula \eqref{PieriRules1}
 \be\label{qWhittakerPieriRule}
P^{qW}_{(m,0,\ldots,0)}(x)\,\times   P^{qW}_\la(x)
  =\,\sum_{\mu_i\geq\la_i\geq\mu_{i+1}\atop |\mu|-|\la|=n}\varphi^q_{\mu/\la}\,
  P^{qW}_\mu(x)\,,\\
  \varphi^q_{\mu/\la}\,
  =\,\Delta_q(\mu)\,\frac{\Theta(\mu_1-\la_1)}{(\mu_1-\la_1)_q!}
  \prod_{i=1}^{\ell}\frac{\Theta(\la_i-\mu_{i+1})}{(\la_i-\mu_{i+1})_q!}
  \frac{\Theta(\mu_{i+1}-\la_{i+1})}{(\mu_{i+1}-\la_{i+1})_q!}\,.
 \ee
The analog of the  Cauchy-Littlewood identity \eqref{CauchyLittlewood}
is given by
 \be\label{qWhittakerCauchyLittlewood}
  \prod_{i=1}^n\prod_{j=1}^m\,\Gamma_q(x_iy_j)\,\,
  =\,\,\sum_{\la \in Y_{n,m}}b^q_{\la}\,
  P^{qW}_{\la}(x)\,P^{qW}_\la(y)\,,\\
  b^q_{\la}\,
  =\,\<P^{qW}_\la,\,P^{qW}_\la\>_q^{-1}\,
  =\,\frac{\Delta_q(\la)}{(\la_{\ell+1})_q!}\,,
 \ee
summed over $\la=(\la_1\geq\ldots\geq\la_m)$ with $m\leq n$.

\begin{rem}
The $q$-Whittaker polynomials $P^{qW}_\la(x)$ can be also obtained
from the Macdonald polynomials $P_{\la}(x)$ under the limit
$t=q^{-k}$, $k\to+\infty$. Let
$$
  D(x)\,=\,\prod_{i=1}^{\ell+1}x_i^{k\varrho_{\ell+2-i}}\,,
$$
$$
  D^{\vee}(\la)\,
  =\,\prod_{i=1}^{\ell+1}q^{\la(\ell\,k+\varrho)}\times
  \prod_{i,j=1\atop i<j}^{\ell+1}
  \frac{1}{\Gamma_{q,\,q^{-k}}(q^{\la_i-\la_j-k(j-i-1)})}
  \frac{\Gamma_{q,\,q^{-k}}(q^{\la_i-\la_j})}
  {\Gamma_{q^{-1},\,q^{-(k+1)}}(q^{\la_j-\la_i})}\,.
$$
Then we have
 \be\label{qWhittakerFunction}
    P^{qW}_{\underline{p}}(q^{\la})\,
    =\,\lim_{k\to+\infty}\Big[\,D^{\vee}(\la)^{-1}\,
    D(q^{\underline{p}+k\varrho})^{-1}\times
    P_{\la+k\rho}(q^{\underline{p}+\varrho(k+1)})\,\Big]\,.
 \ee
for $\underline{p}=(p_1,\ldots ,p_{\ell+1})$ be  a partition
$(p_1\geq\ldots\geq p_{\ell+1})$, and
$P^{qW}_{\underline{p}}(q^{\lambda})=0$ otherwise.
\end{rem}

\subsection{Jack's symmetric polynomials}

Now we consider a bases
of symmetric polynomials consisting  of the Jack polynomials, obtained
from the Macdonald  polynomials by a specialization
(see \cite{M} and references therein).

Let $\kappa$ be a positive integer, and let $\Lambda^{(\ell+1)}$ be
the homogeneous component of degree $\ell+1$ in the algebra of
symmetric functions in variables $x_1,x_2\ldots$ over
$\IQ(\kappa)=\IQ$. Define a pair of scalar products on
$\Lambda^{(\ell+1)}$ using the standard bases $\{p_{\la}\}$ of power
series symmetric polynomials:
 \be
  \<p_{\la},\,p_{\mu}\>_{\kappa}\,
  =\,\kappa^{-l(\la)}\delta_{\la\mu}\,z_{\la}\,,
 \hspace{1.5cm}
  l(\la)\,=\,\bigl|\{m|\,\la_m\neq0\}\bigr|\,,
 \ee
where $z_{\la}=\prod\limits_{n\geq1}n^{m_n}m_n!$ and
$m_n=\bigl|\{k\,:\,\la_k=n\}\bigr|$,
 \be
  \<p_{\la},\,p_{\mu}\>'_{\kappa}\,
  =\,\frac{1}{(\ell+1)!}\int_T\!d^{\times}\!z\,\,
  p_{\la}(z) \,p_{\mu}(z^{-1})\,\Delta_{(\kappa)}(z)\,,
 \ee
where
 \be
  \Delta_{(\kappa)}(z)\,
  =\,\prod_{i,j=1\atop i\neq j}^{\ell+1}
  \bigl(1-z_i^{-1}z_j\bigr)^{\kappa}\,.
 \ee

\begin{de} Jack's symmetric functions $P^{(\kappa)}_{\la}$ are the
elements of $\Lambda^{(\ell+1)}$ such that
$$
 \bigl\<P^{(\kappa)}_{\la},\,P^{(\kappa)}_{\mu}\bigr\>_{\kappa}\,=\,0
$$
whenever $\la\neq\mu$, and
$$
 P^{(\kappa)}_{\la}\,
 =\,m_{\la}\,+\,\sum_{\mu<\la}u_{\la\mu}^{(\kappa)}\,m_{\mu}\,.
$$
\end{de}
The  Jack polynomials $P^{(\kappa)}_{\la}(x)$ are orthogonal with
respect to both scalar products and the following normalization
condition holds:
 \be
  \<P^{(\kappa)}_{\la},\,P^{(\kappa)}_{\la}\>'_{\kappa}\,
  =\,\prod_{i,j=1\atop i<j}^{\ell+1}
  \frac{\Gamma\bigl(\la_i-\la_j+\kappa(j-i+1)\bigr)}
  {\Gamma\bigl(\la_i-\la_j+\kappa(j-i)\bigr)}
  \frac{\Gamma\bigl(\la_i-\la_j+1+\kappa(j-i-1)\bigr)}
  {\Gamma\bigl(\la_i-\la_j+1+\kappa(j-i)\bigr)}\,.
 \ee

Similarly to the cases of the Macdonald polynomials and
$q$-Whittaker polynomials, Jack's polynomials are eigenfunction of
dual families of mutually commuting differential/difference
operators.
\begin{te} {\it (i)}
 The Jack symmetric polynomials are eigenfunctions of a set of mutually
commuting  Sekiguchi differential operators:
 \be\label{JackGenFunctionEigenvalue}
  \CD_{\ell+1}(X)\cdot P^{(\kappa)}_{\la}(x)\,
  =\,\prod_{i=1}^{\ell+1}\bigl(X+(\la_i+\varrho_i\kappa)\bigr)\,
  P^{(\kappa)}_{\la}(x)\,,
 \ee
where
 \be\label{JackGenFunction}
  \CD_{\ell+1}(X)\,=\,\sum_{r=1}^{\ell+1}X^{\ell+1-r}\CH_r\\
  =\,\prod_{i,j=1\atop i<j}^{\ell+1}(x_i-x_j)^{-1}
  \times
  \sum_{\s\in\frak{S}_{\ell+1}}(-1)^{\s}\,
  \prod_{i=1}^{\ell+1}x_i^{\s(\varrho_i)}\Big\{X+\s(\varrho_i)\kappa
  +\,x_i\frac{\pr}{\pr x_i}\Big\}\,,
 \ee
with $\CH_0:=1$.

{\it (ii)} The Jack polynomials
are eigenfunctions of a set of mutually
commuting difference  operators
 \be\label{JackDualGenFunction}
  \CD_{(\kappa)}^{\vee}(X)\cdot P^{(\kappa)}_{\la}(x)\,
  =\,c^{\vee}_{\ell+1}(x;\,X)\,P^{(\kappa)}_{\la}(x)\,,\\
  c^{\vee}_{\ell+1}(x;\,X)\,=\,\prod_{i=1}^{\ell+1}(1+Xx_i)\,,
 \ee
where
$$
 \CD_{(\kappa)}^{\vee}(X)\,=\,\sum_{r=0}^{\ell+1}X^r\CH_r^{\vee},\,
\hspace{1.5cm}
 \CH_0^{\vee}:=1\,,
$$
and
 \be\label{JackDualHamiltonians}
  \CH^{\vee}_r(\la)\,
  =\,\sum_{I_r}\prod_{i\in I_r\atop j\notin I_r}\!
  \,\frac{(i-j+1)\kappa+\la_j-\la_i-1}{(i-j)\kappa+\la_j-\la_i-1}\,
  T^{\vee}_{I_r}\,.
 \ee
Here the summation goes  over all $r$-element subsets $I_r$ of
$(1,2,\ldots,\ell+1)$ and the operators $T^{\vee}$ are defined  as
follows
$$
 T^{\vee}_{I_r}\,=\,\prod_{i\in I_r}T_{q,\,q^{\lambda}}\,,
\hspace{1.5cm}
 T_{q,\,q^{\lambda}}\cdot f_{\lambda}\,
 =\,\left[
 \begin{array}{lc}
 f_{\lambda+\delta_{ij}}\,, & \lambda_{i-1}\geq\lambda_i+1\\
 0 & \mbox{otherwise}
 \end{array}\right.\,,
$$
where
$\lambda+\delta_{ij}\,:=\,\lambda+(0,\ldots,\underbrace{1}_i\,,\ldots,0)$\,.
\end{te}

\begin{rem} The generating function \eqref{JackGenFunction}
can be considered as an appropriately defined non-commutative determinant:
$$
 \CD_{\ell+1}(X)\,
 =\,\frac{1}{\det\|x_i^{\varrho_j}\|}\,
 \det\Big\|x_i^{\varrho_j}\Big(X\,+\,\varrho_j\kappa\,
 +\,x_i\frac{\pr}{\pr x_i}\Big)\Big\|\,.
$$
\end{rem}
In particular, the first statement implies:
 \be\label{JackEigenvalues}
  \CH_r\cdot P^{(\kappa)}_{\la}(x)\,
  =\,\chi_r(\la+\varrho\kappa)\,
  P^{(\kappa)}_{\la}(x)\,,
 \hspace{1.5cm}
  r=1,\ldots,\ell+1\,,
 \ee
where $\varrho_i\,=\,\ell+1-i\,$. The first  two Hamiltonians are
given by
 \be
  \CH_1\,=\,\sum_{i=1}^{\ell+1}\Big\{x_i\frac{\pr}{\pr x_i}\,
  +\,\varrho_i\kappa\Big\}\,,\\
  \CH_2\,=\,\sum_{i,j=1\atop i<j}^{\ell+1}\Big(
  x_i\frac{\pr}{\pr x_i}+\varrho_i\kappa\Big)
  \Big(x_j\frac{\pr}{\pr x_j}+\varrho_j\kappa\Big)\,
  +\,\kappa\sum_{i=1}^{\ell+1}\Big(\varrho_i\,
  +\,\sum_{j=1\atop j\neq i}\frac{x_i}{x_j-x_i}\Big)\,
  x_i\frac{\pr}{\pr x_i}\,.
 \ee
The Jack symmetric functions can be obtained from
Macdonald polynomials by taking the limit $\hbar \to 0$ for
$t=e^{\kappa \h}$, $q=e^{\h}$ \cite{M}
\be\label{JackLimitRelation2}
  \lim_{\hbar\to 0}\Gamma_{q,t}(x)\,
  =\,\frac{1}{(1-x)^{\kappa}}\,
\qquad   \lim_{\hbar \to 0}\,b_{(n)}\,
  =\,\frac{\Gamma(n+\kappa)}{\Gamma(n)\,\Gamma(\kappa)}\,,
\qquad t=e^{\kappa \h},\,\,\,q=e^{\h}.
 \ee
It is easy to infer analogs of \eqref{PieriRules1} and
\eqref{CauchyLittlewood} for Jack polynomials. In particular, the
Pieri rules for the Jack polynomials are given by (see
\cite{S}):
 \be\label{JackPieriRules}
P^{(\kappa)}_{(m)}  \times P^{(\kappa)}_{\la}\,
  =\frac{1}{b^{(\kappa)}_{(m)}}
  \sum_{\mu_i\geq\la_i\geq\mu_{i+1}\atop |\mu|-|\la|=n}
  \varphi^{(\kappa)}_{\mu/\la}\,
  P_{\mu}^{(\kappa)}\,,\\
  \varphi^{(\kappa)}_{\mu/\la}\,\,
  =\,\,\prod_{i,j=1\atop i\leq j}^{\ell+1}\Big[\,
  \frac{\Gamma\bigl(\mu_i-\mu_j+1+(j-i)\kappa\bigr)}
  {\Gamma\bigl(\mu_i-\mu_j+(j-i+1)\kappa\bigr)}
  \frac{\Gamma\bigl(\mu_i-\la_j+(j-i+1)\kappa\bigr)}
  {\Gamma\bigl(\mu_i-\la_j+1+(j-i)\kappa\bigr)}\\
  \hspace{1.5cm}\times\,
  \frac{\Gamma\bigl(\la_i-\la_{j+1}+1+(j-i)\kappa\bigr)}
  {\Gamma\bigl(\la_i-\la_{j+1}+(j-i+1)\kappa\bigr)}
  \frac{\Gamma\bigl(\la_i-\mu_{j+1}+(j-i+1)\kappa\bigr)}
  {\Gamma\bigl(\la_i-\mu_{j+1}+1+(j-i)\kappa\bigr)}\,
  \Big]\,,
 \ee
 and in the product \eqref{JackPieriRules}  we omit the terms
 containing $\lambda_{\ell+2}$ and $\mu_{\ell+2}$.
The analog of the Cauchy-Littlewood identity \eqref{CauchyLittlewood}
is given by
 \be\label{JackCauchyLittlewood}
  \Pi^{(\kappa)}_{n,m}\,:=\,\prod_{i,j}\frac{1}{(1-x_iy_j)^{\kappa}}\,
  =\,\sum_{\la\in Y_{n,m}}\,b^{(\kappa)}_{\la}\,
  P^{(\kappa)}_{\la}(x)\,P^{(\kappa)}_{\la}(y)\,,
 \ee
where  the summation goes over partitions
$\la=(\la_1\geq\ldots\geq\la_{\min(m,n)})$  and
 \be
  b^{(\kappa)}_{\la}\,
  =\,\lim_{\hbar \to 0}b_{\la}\,
  =\,\prod_{(i,j)\in\la}
  \frac{\kappa(\la^\top_j+1-i)+\la_i-j}{\kappa(\la^\top_j-i)+\la_i+1-j}\,,
  \qquad t=e^{\kappa \hbar}, \quad q=e^{\hbar}.
 \ee

\section{Baxter operator formalism for
symmetric polynomials}

In the previous Section we describe various bases in the space of
symmetric polynomials defined as common eigenfunctions of two  sets
of mutually commuting operators called (dual) quantum
Hamiltonians. In this Section we define a dual pair of the  Baxter
operators acting in the space of symmetric polynomials, commuting
with dual pairs of  quantum Hamiltonians. The constructed bases in the space
of polynomials is also a bases of eigenfunctions of the dual pair of
Baxter operators.

\subsection{Baxter operator formalism for
Macdonald symmetric  polynomials}

In this Section we develop the Baxter operator formalism for the
Macdonald polynomials. We construct a dual pair of  Baxter operators
and a dual  pair of recursive operators.  This results in various families
of  integral/sum representations  for the Macdonald  polynomials.

\begin{de} Baxter operator $\CQ_{\gamma}=\CQ_{\gamma}(q,t)$ associated with
Macdonald integrable system is a family of operators  acting on the
space $\Lambda_{q,t}^{(\ell+1)}$ of  symmetric  polynomials as
follows:
 \be\label{BO1}
  \CQ_{\gamma}\cdot P(x)\,
  =\,\int_Td^{\times}\!y\,\,Q_{\gamma}\bigl(x,y\bigr)\,\Delta(y)\,P(y^{-1})\,,
 \hspace{1cm}
  \gamma\in \IZ\,,
 \ee
where integral kernel is given by
 \be
  Q_{\gamma}(x,y)\,
  =\,\prod_{i=1}^{\ell+1}\bigl(x_iy_i\bigr)^{\gamma}
  \prod_{i,j=1}^{\ell+1}\Gamma_{q,t}(x_iy_j)\,.
 \ee
\end{de}

\begin{te}
The Baxter operator \eqref{BO1}
acts on the Macdonald polynomials $P_{\la}(x)$ as follows:
 \be\label{BaxterAction}
  \CQ_{\gamma}\cdot
  P_{\la}(x)\,=\,L_{\gamma}(\la)\,P_{\la}(x)\,,
 \hspace{1.5cm}
  \lambda_{\ell+1}\geq\gamma
\ee
\be
  \CQ_{\gamma}\cdot
  P_{\la}(x)\,=0\,,
 \hspace{3cm}
  \lambda_{\ell+1}<\gamma
\ee
where
\be
  L_{\gamma}(\la)\,:=\,L_{\gamma}(\la,q,t)
  =\,\prod_{i=1}^{\ell+1}
  \frac{\Gamma_{q,\,tq^{-1}}(q)}
  {\Gamma_{q,\,tq^{-1}}(t^{\ell+1-i}q^{\la_i-\gamma+1})}\,.
 \ee
\end{te}
\proof  The  Baxter operator \eqref{BO1}
 can be represented in the following form
 \be\label{BaxterOperator}
  \CQ_{\gamma}\,:=\,D_{\gamma}^{\frak{gl}_{\ell+1}}\circ
  \CC^{\frak{gl}_{\ell+1}}\circ D_{-\gamma}^{\frak{gl}_{\ell+1}}\,,
 \ee
where  the  operator  $\CC^{\frak{gl}_{\ell+1}}$ acts as:
 \be
  \CC^{\frak{gl}_{\ell+1}}\cdot P_{\la}(x)\,:
  =\,\bigl\<\Pi_{\ell+1,\,\ell+1},\,P_{\la}\bigr\>'_{q,\,t}\,,
 \ee
and the operator  $D_{\gamma}^{\frak{gl}_{\ell+1}}$ in
\eqref{BaxterOperator} acts on $P_{\la}(x)$ according the following
rule:
 \be
  D_{\gamma}^{\frak{gl}_{\ell+1}}\cdot P_{\la}(x_1,\ldots,x_{\ell+1})\,:
  =\,P_{\la+(\ell+1)^{\gamma}}(x_1,\ldots,x_{\ell+1})\\
  =\,\,(x_1\cdot\ldots\cdot x_{\ell+1})^{\gamma}\,
  P_{\la}(x_1,\ldots,x_{\ell+1})\,,
 \ee
$P_{\la+(\ell+1)^{\gamma}}(x):= P_{\lambda_1+\gamma,\ldots,\lambda_{\ell+1}+\gamma}(x)$.
The eigenvalue of the  operator $\CC^{\frak{gl}_{\ell+1}}$  can
be found explicitly using the Cauchy-Littlewood identity
\eqref{CauchyLittlewood} and orthogonality of the Macdonald polynomials
with respect to the two scalar products, \eqref{SP0} and
\eqref{AnotherScalarProd}
 \be
  \CC^{\frak{gl}_{\ell+1}}\cdot P_{\la}^{\frak{gl}_{\ell+1}}(x)\,
  =\,\int_Td^{\times}\!y\,
  \Pi_{\ell+1,\ell+1}(x,y)
  P_{\la}^{\frak{gl}_{\ell+1}}(y^{-1})\Delta(y)\\
  =\,\int_Td^{\times}\!y\,
  \Big(\sum_{\mu_1\geq\cdots\geq\mu_{\ell+1}\geq 0}
b_{\mu}P_{\mu}^{\frak{gl}_{\ell+1}}(x)P_{\mu}^{\frak{gl}_{\ell+1}}(y)\Big)
  P_{\la}^{\frak{gl}_{\ell+1}}(y^{-1})\Delta(y)\\

 =\,\sum_{\mu_1\geq\cdots\geq\mu_{\ell+1}\geq 0} \,b_{\mu}\,\Big(\int_Td^{\times}\!y\,
 P_{\mu}^{\frak{gl}_{\ell+1}}(y)
  P_{\la}^{\frak{gl}_{\ell+1}}(y^{-1})\Delta(y)\Big)  P_{\mu}^{\frak{gl}_{\ell+1}}(x)\\
  =\,\sum_{\mu\in\IZ^{\ell+1}}
\Big(\prod_{i=1}^{\ell}\Theta(\mu_i-\mu_{i+1})\Big)\Theta(\mu_{\ell+1})\,
  b_{\mu}\delta_{\lambda,\mu}\,
  \<P_{\mu}^{\frak{gl}_{\ell+1}},\,P_{\la}^{\frak{gl}_{\ell+1}}\>'_{q,\,t}
  P_{\mu}^{\frak{gl}_{\ell+1}}(x)\\
  =\,\prod_{i=1}^{\ell}\Theta(\la_i-\la_{i+1})\Theta(\lambda_{\ell+1})
  \frac{\<P_{\la}^{\frak{gl}_{\ell+1}},\,P_{\la}^{\frak{gl}_{\ell+1}}\>'_{q,\,t}}
  {\<P_{\la}^{\frak{gl}_{\ell+1}},\,P_{\la}^{\frak{gl}_{\ell+1}}\>_{q,\,t}}\,
  P_{\la}^{\frak{gl}_{\ell+1}}(x)\,.
 \ee
Therefore the eigenvalue of the dual Baxter operator
$\CQ_{\gamma}$ on $P_{\lambda}(x)$ is given by
 \be\nn
 \Theta(\la_{\ell+1}-\gamma)
\frac{\<P_{\la-(\ell+1)^{\gamma}},\,P_{\la-(\ell+1)^{\gamma}}\>'_{q,\,t}}
 {\<P_{\la-(\ell+1)^{\gamma}},\,P_{\la-(\ell+1)^{\gamma}}\>_{q,\,t}}\,
 =\,\Theta(\la_{\ell+1}-\gamma)\frac{\<P_{\la},\,P_{\la}\>'_{q,\,t}}
 {\<P_{\la-(\ell+1)^{\gamma}},\,P_{\la-(\ell+1)^{\gamma}}\>_{q,\,t}}\\
 =\,b_{\lambda-(\ell+1)^{\gamma}}
\Theta(\la_{\ell+1}-\gamma)\<P_{\lambda},P_{\lambda}\>'_{q,t}\,.
 \ee
One can rewrite the right hand side of  \eqref{Norm} in the following form:
 \be\label{bLambda}
  b_{\la}\,
  =\,\prod_{i=1}^{\ell+1}
  b_{(\la_i-\la_{i+1})}\,\times\!
  \prod_{i,j=1\atop i<j}^{\ell+1}
  \frac{\Gamma_{q,\,tq^{-1}}(t^{j-i}q^{\la_i-\la_j+1})}
  {\Gamma_{q,\,tq^{-1}}(t^{j-i}q^{\la_i-\la_{j+1}+1})}
  \,,\qquad \la_{\ell+2}:=0\,,
 \ee
where
$$
 b_{(n)}\,=\,\prod_{i=1}^n\frac{1-tq^{n-i}}{1-q^{n+1-i}}\,.
$$
Then combining \eqref{bLambda} with \eqref{AnotherNorm} one readily
arrives at \eqref{BaxterAction}. $\Box$

There is an analog of the classical Baxter equation (see Theorem 2.3
in \cite{GLO1}) relating Baxter operator and quantum Hamiltonian
operators.

\begin{prop} The operators $\CQ_{\gamma}$ given by
\eqref{BO1} and the  generating function $\CM_{\ell+1}(X)$ from
\eqref{MacdonaldGenFunction} commute and satisfy the following
relation:
 \be\label{MacdonaldBaxterEquation}
  \CM_{\ell+1}(-q^{-\gamma})\circ
  \CQ_{\gamma}(q,q^{-k})\,\,
  =\,\,\CQ_{\gamma+1}(q,q^{1-k})\,.
 \ee
\end{prop}
\proof Recall that the Macdonald polynomials are common
eigenfunctions of $\CQ_{\gamma}$ and $\CM_{\ell+1}(X)$. Thus it is
enough to check \eqref{MacdonaldBaxterEquation} on common
eigenvalues of $\CM_{\ell+1}$ and $\CQ_\gamma$ acting on
$P_{\lambda}(x)$. Denoting $L_{\gamma}(\lambda)$ and
$c_{\ell+1}(\lambda)$ the corresponding eigenvalues:
$$
 L_{\gamma}(\la+k\varrho;\,q,\,q^{-k})\, =\,\prod_{i=1}^{\ell+1}\,
 b_{(\la_i-\gamma)}\,,\qquad
 c_{\ell+1}(\lambda+k\varrho;\,-q^{-\gamma})\,
 =\,(1-t)^{-(\ell+1)}\prod_{i=1}^{\ell+1}\bigl(1-q^{\la_i-\gamma}\bigr)\,,
$$
we easily check the following relation:
 \be
  c_{\ell+1}(\la+k\varrho\,;\,-q^{\gamma})\times
  L_{\gamma}(\la+k\varrho\,;\,q,t)\,\,
  =\,\,L_{\gamma+1}(\la+k\varrho\,;\,q,tq)\,.
 \ee
This entails the operator relation
\eqref{MacdonaldBaxterEquation}. $\,\Box$

Now we define the dual  Baxter  operator.

\begin{de} The dual Baxter operator
$\,{}^{^\vee}\!\!\!\!\!\CQ_z=\,{}^{^\vee}\!\!\!\!\!\CQ_z(q,t)$ is a
family of  operators acting in $\Lambda_{q,t}^{(\ell+1)}$ \be
  {}^{^\vee}\!\!\!\!\!\CQ_z\cdot P_{\la}(x)\,
  =\,\sum_{\mu\in\IZ^{\ell+1}}\,
  {}^{^\vee}\!\!\!\!\!Q_z(\la,\,\mu)\,P_{\mu}(x)\,
\ee
with the  kernel function
 \be\label{DualBaxterOperator}
  {}^{^\vee}\!\!\!\!\!Q_z(\la,\,\mu)\,
  =\,z^{|\mu|-|\la|}\,\varphi_{\mu/\la}\,\,,
 \ee
and
 \be
  \varphi_{\mu/\la}\,
  =\,\prod_{i,j=1\atop i\leq j}^{\ell+1}
  \frac{\Gamma_{q,tq^{-1}}\bigl(t^{j-i}q^{\mu_i-\mu_j+1}\bigr)}
  {\Gamma_{q,tq^{-1}}\bigl(t^{j-i}q^{\mu_i-\la_j+1}\bigr)}\,
  \frac{\Gamma_{q,tq^{-1}}\bigl(t^{j-i}q^{\la_i-\la_{j+1}+1}\bigr)}
  {\Gamma_{q,tq^{-1}}\bigl(t^{j-i}q^{\la_i-\mu_{j+1}+1}\bigr)}\\
  \times\Theta(\mu_1-\lambda_1)
  \prod_{i=1}^{\ell}\Theta(\la_i-\mu_{i+1})
  \Theta(\mu_{i+1}-\la_{i+1})\,\,,
 \ee
where in the product one should omit the factors depending on $\la_{\ell+2}$
and $\mu_{\ell+2}$.
\end{de}

\begin{te} The action of the dual Baxter operator on the Macdonald
polynomials reads
 \be\label{DualBaxterAction}
  {}^{^\vee}\!\!\!\!\!\CQ_z\cdot P_{\la}(x)\,
  =\,L^{\vee}_z(x)\,P_{\la}(x)\,,
\ee
where the eigenvalue is given by
\be
  L^{\vee}_z(x)\,
  =\,\prod_{i=1}^{\ell+1}\Gamma_{q,t}(zx_i)\,.
 \ee
\end{te}
\proof The statement of the Theorem directly follows from the Pieri
formula for the Macdonald polynomials \eqref{PieriRules1},
\eqref{PieriRules2}:
 \be
  {}^{^\vee}\!\!\!\!\!\CQ_z\cdot P_{\la}(x)\,
  =\,\sum_{\mu_i\geq\la_i\geq\mu_{i+1}}\,z^{|\mu|-|\la|}\,\varphi_{\mu/\la}\,
  P_{\mu}(x)\\
  =
  \,\sum_{m=0}^{\infty}z^m\sum_{{\mu_i\geq\la_i\geq\mu_{i+1}}\atop |\mu|-|\la|=m}
  \varphi_{\mu/\la}\,P_{\mu}(x)\,=\,
  \,\Big(\sum_{m=0}^{\infty} b_{(m)}z^m P_m(x)\Big)P_{\lambda}(x)\\

  =\,\Pi_{\ell+1,\,1}(x,z)\,P_{\la}(x)\,,
 \ee
and identification  $L^{\vee}_z(x)=\Pi_{\ell+1,\,1}(x,z)$.  $\Box$

\begin{prop} The operator $\,{}^{^\vee}\!\!\!\!\!\CQ_z$
satisfies the following difference relation:
 \be\label{MacdonaldDualBaxterEquation}
  \CM^{\vee}_{\ell+1}(-q^{k\ell/2}z)\circ\,{}^{^\vee}\!\!\!\!\!\CQ_z(q,q^{-k})\,\,
  =\,\,\,{}^{^\vee}\!\!\!\!\!\CQ_{qz}\bigl(q,q^{-k-1}\bigr)\,.
 \ee
where $\CM^{\vee}_{\ell+1}(X)$ is the  generating function
\eqref{MacdonaldDualGenFunction} of the dual quantum Hamiltonians.

\end{prop}
\proof It is enough to check \eqref{MacdonaldDualBaxterEquation} on
common eigenfunctions $P_{\lambda}(x)$ of $\CM^{\vee}_{\ell+1}$ and
$\CQ^{\vee}$.  Denoting  $L^{\vee}_{\mu}$ and $c^{\vee}_{\ell+1}$
the corresponding eigenvalues
$$
 L_z^{\vee}(x)\,=\,\prod_{i=1}^{\ell+1}\Gamma_{q,t}(zx_i)\,,\qquad
 c^{\vee}_{\ell+1}(x;\,-t^{-\ell/2}z)\,=\,\prod_{i=1}^{\ell+1}(1-zx_i)\,,
$$
we easily check the following relation:
 \be
  c^{\vee}_{\ell+1}(x;\,-t^{-\ell/2}z)\times
  L^{\vee}_{z}(x\,;\,q,t)\,\,
  =\,\,L^{\vee}_{qz}(x\,;\,q,tq^{-1})\,.
 \ee
This entails the operator relation \eqref{MacdonaldDualBaxterEquation}. $\,\Box$

Let us introduce the following notation
$$
 \underline{a}_n\,=\,(a_{n,1},\ldots,a_{n,\,n})\,,
\hspace{1.5cm}
 \underline{a}_n'\,=\,(a_{n,1},\ldots,a_{n,\,n-1}).
$$
The following recursive relations  hold; the first one (see
\cite{AOS})
 \be\label{MacdonaldRecursion}
  P_{\la}(\underline{x}_{\ell+1})\,
  =\,\int_T\!d^{\times}\!\underline{x}_{\ell}\,
  Q^{\frak{gl}_{\ell+1}}_{\frak{gl}_{\ell}}(\underline{x}_{\ell+1};\,
  \underline{x}_{\ell}|\,\la_{\ell+1})\,
  \Delta(\underline{x}_{\ell})\,P_{\la'}(\underline{x}_{\ell}^{-1})\,,\\
  Q^{\frak{gl}_{\ell+1}}_{\frak{gl}_{\ell}}(\underline{x}_{\ell+1};\,
  \underline{x}_{\ell}|\,\la_{\ell+1})\,
  =\,x_{\ell+1,\,\ell+1}^{\la_{\ell+1}}\prod_{i=1}^{\ell}
  \bigl(x_{\ell+1,\,i}x_{\ell,\,i}\bigr)^{\la_{\ell+1}}\times
  \Pi_{\ell+1,\,\ell}(x_{\ell+1,\,i},\,x_{\ell,\,i})\,,
 \ee
 and the dual recursive relation (see \cite{M}):
 \be\label{MacdonaldDualRecursion}
  P_{\underline{\la}_{\ell+1}}(x)\,
  =\,\sum_{\underline{\la}_{\ell}}\,{}^{^\vee}\!\!\!\!\!
  Q^{\frak{gl}_{\ell+1}}_{\frak{gl}_{\ell}}(\underline{\la}_{\ell+1},\,
  \underline{\la}_{\ell}|\,x_{\ell+1})\,
  P_{\underline{\la}_{\ell}}(x')\,,\\
  Q^{\frak{gl}_{\ell+1}}_{\frak{gl}_{\ell}}(\underline{\la}_{\ell+1},\,
  \underline{\la}_{\ell}|\,x_{\ell+1})\,
  =\,x_{\ell+1}^{|\underline{\la}_{\ell+1}|-|\underline{\la}_{\ell}|}
  \psi_{\underline{\la}_{\ell+1}/\underline{\la}_{\ell}}\,,
 \ee
where
 \be
  \psi_{\la/\mu}\,
  =\,\prod_{1\leq i\leq j\leq\ell}
  \frac{\Gamma_{q,\,tq^{-1}}(t^{j-i}q^{\mu_i-\mu_j+1})}
  {\Gamma_{q,\,tq^{-1}}(t^{j-i}q^{\la_i-\mu_j+1})}
  \frac{\Gamma_{q,\,tq^{-1}}(t^{j-i}q^{\la_i-\la_{j+1}+1})}
  {\Gamma_{q,\,tq^{-1}}(t^{j-i}q^{\mu_i-\la_{j+1}+1})}
 \ee
when $\la$ and $\mu$ are interlaced (i.e.
$\la_1\geq\mu_1\geq\ldots\geq\la_{\ell}\geq\mu_{\ell}\geq\la_{\ell+1}\geq
0$), and $\psi_{\la/\mu}=0$ otherwise.

These  recursive relations allow to introduce the corresponding
recursive operators $\CQ^{\frak{gl}_{n+1}}_{\frak{gl}_n}(\la_{n+1})$
and
$\,{}^{^\vee}\!\!\!\!\!\CQ^{\frak{gl}_{n+1}}_{\frak{gl}_n}(x_{n+1})$:
 \be
\CQ^{\frak{gl}_{\ell+1}}_{\frak{gl}_{\ell}}(\la_{\ell+1})\cdot
f(\underline{x}_{\ell+1})\,
  =\,\int_T\!d^{\times}\!\underline{x}_{\ell}\,
  Q^{\frak{gl}_{\ell+1}}_{\frak{gl}_{\ell}}(\underline{x}_{\ell+1};\,
  \underline{x}_{\ell}|\,\la_{\ell+1})\,
  \Delta(\underline{x}_{\ell})\,f(\underline{x}_{\ell}^{-1})\,,
 \ee
and
 \be
  {}^{^\vee}\!\!\!\!\!
  \CQ^{\frak{gl}_{\ell+1}}_{\frak{gl}_{\ell}}(x_{\ell+1})\cdot
 f(\underline{\la}_{\ell+1})\,
  =\,\sum_{\underline{\la}_{\ell}}\,{}^{^\vee}\!\!\!\!\!
  Q^{\frak{gl}_{\ell+1}}_{\frak{gl}_{\ell}}(\underline{\la}_{\ell+1},\,
  \underline{\la}_{\ell}|\,x_{\ell+1})\,
 f(\underline{\la}_{\ell})\,,
 \ee
The existence of the  two dual recursive
representations, \eqref{MacdonaldRecursion} and
\eqref{MacdonaldDualRecursion}, provide a family of $2^{\ell}$
integral representations for the Macdonald polynomials. Namely, let
us change our notations as follows:
$$
 R_{n+1,\,n}^{I}\,:
 =\,\CQ^{\frak{gl}_{n+1}}_{\frak{gl}_n}(\la_{n+1})\,,
\hspace{1cm}
 R_{n+1,\,n}^{II}\,:
 =\,\,{}^{^\vee}\!\!\!\!\!\CQ^{\frak{gl}_{n+1}}_{\frak{gl}_n}(x_{n+1})\,,
\hspace{1.5cm}
 n=1,\ldots,\ell\,;
$$
then for every array $\epsilon=(\epsilon_1,\ldots,\epsilon_{\ell})$
of $\epsilon_n\in\{I,\,II\},\,n=1,\ldots,\ell$ the following
holds:
 \be\label{MixedReps}
  P_{\underline{\la}_{\ell+1}}^{\frak{gl}_{\ell+1}}(\underline{x}_{\ell+1})\,
  =\,\Big\{\,R_{\ell+1,\,\ell}^{\epsilon_\ell}\circ\ldots\circ
  R_{2,\,1}^{\epsilon_1}\circ
  R_{1,\,0}\Big\}\cdot1\,,
 \ee
where $R_{1,\,0}$ is the ${\frak{gl}_1}$-Macdonald polynomial
$P^{\frak{gl}_1}$.

Let us remark that the recursive operators can be factorized
into Baxter operators, similarly to recursive operators for
$\frak{gl}_{\ell+1}$-Whittaker function (see Proposition 3.3 in
\cite{GLO1}). This reveals the fundamental role of the Baxter
operators in the description of various bases of symmetric
polynomials.

\subsection{Baxter operator formalism for
class one  $q$-Whittaker functions}

Now we provide similar results for  $q$-Whittaker polynomials
$P^{qW}_{\la}(x)$.
\begin{de} The Baxter operator acting in $\Lambda_q^{(\ell+1)}$ is
a family of  integral operators
 \be\label{qWhitBaxterOperator}
  {\bf Q}_z\cdot f(\la)\,
  =\,\sum_{\mu\in\IZ^{\ell+1}}\,\Delta_q(\mu)\,
  {\bf Q}_{\ell+1,\,\ell+1}(\mu;\,\la|\,z)\,
  f(\mu)\,,
\ee
with the kernel
 \be\label{qWhhitBaxterKernel}
  {\bf Q}_{\ell+1,\,\ell+1}(\mu,\,\la|\,z)\,
  =\,z^{|\mu-\la|}\,\varphi^q_{\mu/\la}\,
  =\,z^{|\mu-\la|}\,\varphi_{\mu/\la}(q,\,t=0)\times
  \Delta_q(\la)^{-1}\\
  =\,z^{|\mu-\la|}\,
  \,\frac{\Theta(\mu_1-\la_1)}{(\mu_1-\la_1)_q!}
  \prod_{i=1}^{\ell}\frac{\Theta(\la_i-\mu_{i+1})}{(\la_i-\mu_{i+1})_q!}
  \frac{\Theta(\mu_{i+1}-\la_{i+1})}{(\mu_{i+1}-\la_{i+1})_q!}\,.
 \ee
\end{de}

\begin{te} {\it (i)} The action of the Baxter operator ${\bf Q}_z$ on
$q$-Whittaker polynomials \eqref{qWhittakerFunction} is given by
 \be\label{qWhittakerBaxterAction}
  {\bf Q}_z\cdot P^{qW}_{\la}(x)\,
  =\,{\bf L}_z(x)\,P^{qW}_{\la}(x)\,,
 \ee
where
 \be
  {\bf L}_z(x)\,=\,\prod_{i=1}^{\ell+1}\Gamma_q(zx_i)\,.
 \ee
{\it (ii)} The operators ${\bf Q}_z$ and $D^q_{\ell+1}(X)$
\eqref{qWhitEigenvalues}  satisfy the following relation:
 \be\label{qWhittakerBaxterEquation}
  D_{\ell+1}(-z)\circ{\bf Q}_z\,\,
  =\,\,{\bf Q}_{qz}\,.
 \ee
\end{te}
\proof The relation \eqref{qWhittakerBaxterAction} is a direct
consequence of the Pieri formula \eqref{qWhittakerPieriRule}. The
relation \eqref{qWhittakerBaxterEquation} follows from the  relation
between the corresponding eigenvalues: ${\bf
c}^q_{\ell+1}(x;\,-z)\,{\bf L}_z(x)={\bf L}_{qz}(x)$. $\Box$

\begin{de} The dual Baxter operator acting in $\Lambda_q^{(\ell+1)}$
is a family of  integral operators
 \be\label{qWhitDualBaxterOperator}
  {}^{^{\vee}}\!\!\!\!\!{\bf Q}_{\g}\cdot P(x)\,
  =\,\int_Td^\times y\,\,
  {\bf Q}^{\vee}_{\ell+1,\,\ell+1}(x,\,y;\,\g)\,
  \Delta^{\vee}_q(y)\,\,P(y^{-1})\,,
 \hspace{1cm}
  \gamma\in\IZ\,,
 \ee
with the kernel
 \be\label{qWhhitDualBaxterKernel}
  {}^{^{\vee}}\!\!\!\!\!{\bf Q}_{\ell+1,\,\ell+1}(x,\,y;\,\g)\,
  =\,\prod_{i,j=1}^{\ell+1}\bigl(x_iy_i)^{\g}
  \Gamma_q(x_iy_j)\,.
 \ee
\end{de}

\begin{te} {\it (i)} The action of the dual Baxter operator $\,{}^{^{\vee}}\!\!\!\!\!{\bf
Q}_{\g}$ on $q$-Whittaker polynomials reads as follows:
 \be
  {}^{^{\vee}}\!\!\!\!\!{\bf Q}_{\g}\cdot P^{qW}_{\la}(x)\,
    =\,\,{\bf L}^{\vee}_{\g}(\la)\,
  P^{qW}_{\la}(x)\,,
 \hspace{1.5cm}
  {\bf L}^{\vee}_{\g}(\la_1,\ldots,\la_{\ell+1})\,
  =\,\frac{1}{(\la_{\ell+1}-\g)_q!}\,,
 \ee
when $\g\leq\la_{\ell+1}$, and
 \be
  {}^{^{\vee}}\!\!\!\!\!{\bf Q}_{\g}\cdot P^{qW}_{\la}(x)\,
  =\,0\,,
 \hspace{1.5cm}
  \g>\la_{\ell+1}\,.
 \ee
{\it (ii)} The dual Baxter operator satisfies
the following difference equation:
 \be\label{qWhittakerDualBaxterEquation}
  \Big\{1\,-\,{q^{-\g}}\,{\bf H}_1^{\vee}\Big\}\circ
  \,{}^{^{\vee}}\!\!\!\!\!{\bf Q}_{\g}\,
  =\,\,{}^{^{\vee}}\!\!\!\!\!{\bf Q}_{\g+1}\,\,.
 \ee
\end{te}
\proof The first statement follows from
\eqref{qWhittakerCauchyLittlewood} and the orthogonality of the
$q$-Whittaker polynomials (see \cite{GLO1}). The second statement
follows from the relation between the corresponding eigenvalues:
$(1-q^{\la_{\ell+1}-\g})\times
 {\bf L}^{\vee}_{\g}(\la)
 ={\bf L}^{\vee}_{\g+1}(\la)\,.$ $\Box$

In \cite{GLO2} and \cite{GLO4} the following recursive relations for
the $q$-deformed Whittaker functions were established.

\begin{prop}[\cite{GLO2},\cite{GLO4}] The following recursive
relations hold:
 \be\label{qWhittakerRecursion}
  P^{qW}_{\la}(x)\,
  =\,\bigl({\bf Q}^{\frak{gl}_{\ell+1}}_{\frak{gl}_{\ell}}(x_{\ell+1})
  \cdot
  P^{qW}(x')\bigr)_{\la}\\
  =\,\sum_{\la_i\geq\mu_i\geq\la_{i+1}}\!\!\!
  {\bf Q}_{\ell+1,\,\ell}(\la;\,\mu|\,x_{\ell+1})\,
  \Delta_q(\mu)\,
  P^{qW}_{\mu}(x)\,,
 \ee
where
 \be
  {\bf Q}_{\ell+1,\,\ell}(\la;\,\mu|\,z)\,
  =\,z^{|\la|-|\mu|}\prod_{i=1}^{\ell}
  \frac{\Theta(\la_{i}-\mu_i)}{(\la_{i}-\mu_i)_q!}
  \frac{\Theta(\mu_i-\la_{i+1})}{(\mu_i-\la_{i+1})_q!}\,.
 \ee
and
 \be\label{qWhittakerDualRecursion}
  P^{qW}_{\la}(x)\,
  =\,\,{}^{^\vee}\!\!\!\!\!{\bf Q}^{\frak{gl}_{\ell+1}}_{\frak{gl}_{\ell}}(\la_{\ell+1})
  \cdot
  P^{qW}_{\la'}(x)\,
  =\,\int_T\!d^{\times}\!y\,
  \,{}^{^{\vee}}\!\!\!\!\!{\bf Q}_{\ell+1,\,\ell}(x;\,
  y|\,\la_{\ell+1})\,\Delta^{\vee}_q(y)\,
  P^{qW}_{\la'}(y^{-1})\,.
 \ee
where
 \be
  {}^{^{\vee}}\!\!\!\!\!{\bf Q}_{\ell+1,\,\ell}(x_1,\ldots,x_{\ell+1};\,
  y_1,\ldots,y_{\ell}|\,\g)\,
  =\,\prod_{i=1}^{\ell+1}x_i^{\g}\prod_{j=1}^{\ell}y_j^{\g}
  \prod_{1\leq i\leq\ell+1\atop 1\leq j\leq\ell}\Gamma_q(x_iy_j)\,.
 \ee
\end{prop}

The action of recursive operators \eqref{qWhittakerRecursion} and
\eqref{qWhittakerDualRecursion} provide a pair of dual integral/sum
representations of the $q$-deformed Whittaker functions (see
\cite{GLO2} and \cite{GLO4}). Combining the recursive operators of
different types one can obtain $2^{\ell}$  explicit formulas for
$q$-Whittaker functions similarly to \eqref{MixedReps}.

\begin{rem} The recursive operators ${\bf
Q}^{\frak{gl}_{\ell+1}}_{\frak{gl}_{\ell}}(x_{\ell+1})$ and
$\,{}^{^\vee}\!\!\!\!\! {\bf
Q}^{\frak{gl}_{\ell+1}}_{\frak{gl}_{\ell}}(\la_{\ell+1})$ can be
factorized into the Baxter operators \eqref{qWhitBaxterOperator} and
\eqref{qWhitDualBaxterOperator}, similarly to Proposition 3.3 from
\cite{GLO1}.
\end{rem}

\subsection{Baxter operator formalism for Jack's symmetric polynomials}

Now we consider Baxter
operator formalism associated with the Jack symmetric polynomials.

\begin{de}  Baxter operator
is a family of integral operators  acting in
$\Lambda^{(\ell+1)}$ by
 \be\label{JackBaxterOperator}
  \CQ_{\g}^{(\kappa)}\cdot P^{(\kappa)}_{\la}(x)\,
  =\,\int_Td^{\times}y\,\,Q^{(\kappa)}(x,y;\,\g)\,
  \Delta_{(\kappa)}(y)\,P^{(\kappa)}_{\la}(y^{-1})\,,
 \hspace{1cm}
  \g\in\IZ\,,
\ee
with the kernel
\be\label{JackBaxterKernel}
  Q^{(\kappa)}_{\g}(x,y;\,\g)\,
  =\,\prod_{i=1}^{\ell+1}(x_iy_i)^{\g}\,
  \Pi^{(\kappa)}_{\ell+1,\,\ell+1}(x,y)\,.
 \ee
\end{de}

\begin{te} {\it (i)} The action of the Baxter operator  $\CQ_{\g}^{(\kappa)}$ on
the Jack polynomials is given by
 \be\label{JackBaxterAction}
  \CQ_{\g}^{(\kappa)}\cdot P^{(\kappa)}_{\la}(y)\,
  =\,\CL_{\g}(\la)\,P^{(\kappa)}_{\la}(x)\,,
 \hspace{1.5cm}
  \CL_{\g}(\la)\,
  =\,\prod_{i=1}^{\ell+1}
  \frac{\Gamma(\la_i-\g+(\varrho_i+1)\kappa)}
  {\Gamma(\la_i-\g+\varrho_i\kappa+1)}\,,
 \ee
when $\g\leq\la_{\ell+1}+\kappa$, with
$\varrho_i=\ell+1-i\,,i=1,\ldots,\ell+1$, and
 \be
  \CQ_{\g}^{(\kappa)}\cdot P^{(\kappa)}_{\la}(y)\,
  =\,0,\,
 \hspace{1.5cm}
  \g>\la_{\ell+1}+\kappa\,.
 \ee

{\it (ii)} The Baxter operator
$\CQ_{\g}^{(\kappa)}$ commutes with the generating function
$\CD_{\ell+1}$ \eqref{JackGenFunction} and satisfies the following
difference equation:
 \be\label{JackBaxterEquation}
  \CD_{\ell+1}(\kappa-\g)\circ
  \CD_{\ell+1}(1-\g)^{-1}\circ\CQ_{\g}^{(\kappa)}\,
  =\,\CQ_{\g-1}^{(\kappa)}\,.
 \ee
\end{te}
\proof \eqref{JackBaxterAction} follows from
\eqref{JackCauchyLittlewood} and the orthogonality of the Jack
polynomials \cite{M}. $\Box$

\begin{de} The dual Baxter operator $\,{}^{^\vee}\!\!\!\!\!\CQ_z^{(\kappa)}$
is a family of operators acting in $\Lambda^{(\ell+1)}$ by
 \be\label{JackDualBaxterOperator}
  \,{}^{^\vee}\!\!\!\!\!\CQ_z^{(\kappa)}\cdot P^{(\kappa)}_{\la}(x)\,
  =\,\sum_{\mu\in\IZ^{\ell+1}}\,
  {}^{^\vee}\!\!\!\!\!Q^{(\kappa)}(\mu,\la;\,z)\,
  P^{(\kappa)}_{\mu}(x),
 \ee
with the kernel
 \be\label{JackDualBaxterKernel}
  {}^{^\vee}\!\!\!\!\!Q^{(\kappa)}(\mu,\la;\,z)\,
  =\,z^{|\mu|-|\la|}\,\prod_{i,j=1\atop i\leq j}^{\ell+1}\Big[\,
  \frac{\Gamma\bigl(\mu_i-\mu_j+1+(j-i)\kappa\bigr)}
  {\Gamma\bigl(\mu_i-\mu_j+(j-i+1)\kappa\bigr)}
  \frac{\Gamma\bigl(\mu_i-\la_j+(j-i+1)\kappa\bigr)}
  {\Gamma\bigl(\mu_i-\la_j+1+(j-i)\kappa\bigr)}\\
  \hspace{1.5cm}\times\,
  \frac{\Gamma\bigl(\la_i-\la_{j+1}+1+(j-i)\kappa\bigr)}
  {\Gamma\bigl(\la_i-\la_{j+1}+(j-i+1)\kappa\bigr)}
  \frac{\Gamma\bigl(\la_i-\mu_{j+1}+(j-i+1)\kappa\bigr)}
  {\Gamma\bigl(\la_i-\mu_{j+1}+1+(j-i)\kappa\bigr)}\,
  \Big]\\
  \times\Theta(\mu_1-\la_1)\prod_{i=1}^{\ell}\Theta(\la_i-\mu_{i+1})\,
  \Theta(\mu_{i+1}-\la_{i+1})\,.
 \ee
\end{de}

\begin{te} {\it (i)} The action of the dual Baxter operator  on
the Jack polynomials is given by
 \be\label{JackDualBaxterAction}
  \,{}^{^\vee}\!\!\!\!\!\CQ_z^{(\a)}\cdot P^{(\kappa)}_{\la}(x)\,
  =\,\,\,{}^{^\vee}\!\!\!\!\!\CL_z(x)\,P^{(\kappa)}_{\la}(x)\,,
 \ee
where
 \be\label{JackDualBaxterEquation}
  \,\,{}^{^\vee}\!\!\!\!\!\CL_z(x)\,
  =\,\prod_{i=1}^{\ell+1}\frac{1}{(1-zx_i)^{\kappa}}\,.
 \ee
{\it (ii)} The dual Baxter operator commutes
with $ \CD_{(\kappa)}^{\vee}$ given by \eqref{JackDualGenFunction}
and satisfies the following difference equation:
 \be
  \CD_{(\kappa)}^{\vee}(-z)\circ\,{}^{^\vee}\!\!\!\!\!\CQ_z^{(\kappa)}\,\,
  =\,\,\,{}^{^\vee}\!\!\!\!\!\CQ_z^{(\kappa-1)}\,.
 \ee
\end{te}
\proof The first relation   \eqref{JackDualBaxterAction}
follows from the Pieri formula \eqref{JackPieriRules}, and the
Cauchy-Littlewood identity \eqref{JackCauchyLittlewood}  is implied
by the relation between the eigenvalues:
$c^{\vee}_{\ell+1}(x;\,-z)\times\,\,{}^{^\vee}\!\!\!\!\!\CL^{(\kappa)}_z(x)\,
=\,\,\,{}^{^\vee}\!\!\!\!\!\CL^{(\kappa-1)}_z(x)$. $\Box$

The recursive relations  \eqref{MacdonaldRecursion} and
\eqref{MacdonaldDualRecursion} imply similar recursive relations for
Jack's symmetric functions.  Namely, the following recursive
relation hold (see \cite{AMOS1}, \cite{AMOS2}, \cite{AOS}):
 \be\label{JackIntegralRepresentation}
  P_{\underline{\la}_{\ell+1}}^{(\kappa)}(\underline{x}_{\ell+1})\,
  =\,\CQ^{\frak{gl}_{\ell+1}}_{\frak{gl}_{\ell}}(\la_{\ell+1,\,\ell+1})
  \cdot
  P_{\underline{\la}_{\ell+1}'}^{(\kappa)}\\
  =\,\int_T\!d^{\times}x_{\ell}\,\,
  Q^{(\kappa)}_{\ell+1,\,\ell}(\underline{x}_{\ell+1};\,\underline{x}_{\ell}|\,
  \la_{\ell+1,\,\ell+1})\,
  \Delta_{(\kappa)}(\underline{x}_{\ell})\,
  P_{\underline{\la}_{\ell+1}'}^{(\kappa)}(\underline{x}_{\ell}^{-1}),
 \ee
where
 \be\label{JackRecursiveKernel}
  Q^{(\kappa)}_{\ell+1,\,\ell}(\underline{x}_{\ell+1};\,\underline{x}_{\ell}|\,
  \la_{\ell+1,\,\ell+1})\\
  =\,x_{\ell+1,\,\ell+1}^{\la_{\ell+1,\,\ell+1}} \prod_{1\leq i\leq\ell}
  (x_{\ell+1,\,i}x_{\ell,\,i})^{\la_{\ell+1,\,\ell+1}}\times
  \prod_{i=1}^{\ell+1}\prod_{j=1}^{\ell}
  \frac{1}{\bigl(1-x_{\ell+1,\,i}x_{\ell,\,j}\bigr)^{\kappa}}\,.
 \ee

We also have the dual recursive relations
 \be\label{JackDualIntegralRepresentation}
  P^{(\kappa)}_{\underline{\la}_{\ell+1}}(\underline{x}_{\ell+1})\,
  =\,\bigl(\,{}^{^{\vee}}\!\!\!\!
  \CQ^{\frak{gl}_{\ell+1}}_{\frak{gl}_{\ell}}(x_{\ell+1,\,\ell+1})
  \cdot
  P^{(\kappa)}(\underline{x}_{\ell+1}')\,\bigr)_{\underline{\la}_{\ell}}\\
  =\,\sum_{\la_{\ell+1,\,i}\geq\la_{\ell,\,i}\geq\la_{\ell+1,\,i+1}}
  {}^{^{\vee}}\!\!\!\!Q^{\frak{gl}_{\ell+1}}_{\frak{gl}_{\ell}}\bigl(
  \underline{\la}_{\ell+1};\,\underline{\la}_{\ell}|\,x_{\ell+1,\,\ell+1}\bigr)\,
  P^{(\kappa)}_{\underline{\la}_{\ell}}(\underline{x}_{\ell+1}').
 \ee
Here
 \be\label{JackDualRecursiveKernel}
  {}^{^\vee}\!\!\!\!\!Q^{\frak{gl}_{\ell+1}}_{\frak{gl}_{\ell}}\bigl(
  \mu;\,\la|\,
  z\bigr)\\
  =\,z^{|\mu|-|\la|}\,
  \prod_{i,j=1\atop i\leq j}^{\ell}
  \Big[\,
  \frac{\Gamma\bigl(\mu_i-\mu_j+1+(j-i)\kappa\bigr)}
  {\Gamma\bigl(\mu_i-\mu_j+(j-i+1)\kappa\bigr)}
  \frac{\Gamma\bigl(\mu_i-\la_j+(j-i+1)\kappa\bigr)}
  {\Gamma\bigl(\mu_i-\la_j+1+(j-i)\kappa\bigr)}\\
  \hspace{1.5cm}\times\,
  \frac{\Gamma\bigl(\la_i-\la_{j+1}+1+(j-i)\kappa\bigr)}
  {\Gamma\bigl(\la_i-\la_{j+1}+(j-i+1)\kappa\bigr)}
  \frac{\Gamma\bigl(\la_i-\mu_{j+1}+(j-i+1)\kappa\bigr)}
  {\Gamma\bigl(\la_i-\mu_{j+1}+1+(j-i)\kappa\bigr)}\,
  \Big]\,,
 \ee
when $(\mu_1,\ldots,\mu_{\ell+1})$ and $(\la_1,\ldots,\la_{\ell},0)$
are interlaced, and
$\,{}^{^\vee}\!\!\!\!\!Q^{\frak{gl}_{\ell+1}}_{\frak{gl}_{\ell}}(
  \mu;\,\la|\,z)=0$ otherwise.
Obviously for Jack polynomials one has the proper analogs of the
mixed integral/sum representations \eqref{MixedReps} of the
Macdonald polynomials.

\begin{ex} The simplest dual recursive operator  intertwining
$\frak{gl}_2$ and $\frak{gl}_1$ Jack's symmetric functions reads as
follows:
 \be
  {}^{^\vee}\!\!\!\!\!Q^{\frak{gl}_2}_{\frak{gl}_1}
  (\la_{21},\,\la_{22};\,\la_{11}\,|\,x_2)\\
  =\,x_2^{\la_{21}+\la_{22}-\la_{11}}\,
  \frac{\Gamma(\kappa+\la_{21}-\la_{11})\,
  \Gamma(\kappa+\la_{11}-\la_{22})}{\Gamma(\kappa+\la_{21}-\la_{22})}
  \frac{(\la_{21}-\la_{22})!}{(\la_{21}-\la_{11})!\,(\la_{11}-\la_{22})!}\,.
 \ee
This leads to the following
representation of the  $\frak{gl}_2$-Jack's polynomial:
$$
 P^{(\kappa)}_{\la_1,\,\la_2}(x_1,\,x_2)\,
 =\,\sum_{\mu=\la_2}^{\la_1}
 \frac{\Gamma(\kappa+\la_1-\mu)\,
 \Gamma(\kappa+\mu-\la_2)}{\Gamma(\kappa+\la_1-\la_2)}
 \frac{(\la_1-\la_2)!}{(\la_1-\mu)!\,(\mu-\la_2)!}\,
 x_1^{\mu}x_2^{\la_1+\la_2-\mu}\,.
$$
\end{ex}

\begin{rem} The recursive operators
$\CQ^{\frak{gl}_{\ell+1}}_{\frak{gl}_{\ell}}(\la_{\ell+1,\,\ell+1})$
and
$\,{}^{^{\vee}}\!\!\!\!\CQ^{\frak{gl}_{\ell+1}}_{\frak{gl}_{\ell}}(x_{\ell+1,\,\ell+1})$
can be factorized into the Baxter operators \eqref{JackBaxterOperator},
\eqref{JackDualBaxterOperator}, similarly to Proposition 3.3 from
\cite{GLO1}.
\end{rem}

\section{Appendix:  Various analogs of classical $\Gamma$-function}

In this Appendix we provide basic facts on the analogs of classical
$\Gamma$-function arising in the Baxter operator formalism for
Macdonald, $q$-Whittaker  and Jack polynomials.

 Classical $\Gamma$-function can be defined by analytic
continuation of the function defined by Euler's integral
representation:
 \be\label{GammaEulerRep}
  \Gamma(s)\,=\,\int_{\IR}\!\!dt\,\,e^{st}e^{-e^t}\,,
 \hspace{1.5cm}
  {\rm Re}(s)>0.
 \ee
Equivalently $\Gamma$-function is defined as a solution of the
functional equation
$$
\Gamma(s+1)=s\Gamma(s),\qquad \Gamma(1)=1,
$$
such that $\frac{1}{\Gamma(s)}$ is an entire function on the  complex
plane. $\Gamma$-function  allows a representation as the
Weierstrass  product
 \be\label{GammaWeierstrass}
  \Gamma(1+s)\,
  =\,e^{-\gamma s}\prod_{n=1}^{\infty}
  e^{\frac{s}{n}}\,\Big(1+\frac{s}{n}\Big)^{-1}\,,
 \ee
where $\gamma=-\Gamma'(1)$ is the Euler constant. Also the following
reflection property holds
 \be\label{GammaFunctionalEquation}
  \Gamma(s)\times\Gamma(1-s)\, =\,\frac{\pi}{\sin(\pi s)}\,.
 \ee
Note that the integral representation  \eqref{GammaEulerRep} can be inverted
via the Mellin transform
 \be\label{InversI}
  e^{-e^\tau}\,
  =\,\frac{1}{2\pi \imath}\,\int\limits_{\imath\IR+\e}\!\!ds\,\,
  e^{-s\tau}\Gamma(s)\,.
 \ee

Define $(q,t)$-analog of the classical $\Gamma$-function as the
following infinite product
 \be\label{qtG}
  \Gamma_{q,t}(x)
  =\,\prod_{n=0}^{\infty}\frac{1-txq^n}{1-xq^n}\,,
 \ee
where we imply that $q$ is variable taking values in $|q|<1$. This
function has poles at $x=q^{-m}$, $m\in \IZ_{\geq 0}$ and zeroes at
$x=t^{-1}q^{-m}$, $m\in \IZ_{\geq
  0}$.

The function  $\Gamma_{q,t}$ defined by \eqref{qtG} possesses all
the basic properties of the classical $\Gamma$-function outlined
above. The analog of the Weierstrass product is given by
\eqref{qtG}. The analog of the relation
\eqref{GammaFunctionalEquation} is given by:
 \be\label{qtReflectionEquation}
  \Gamma_{q,\,t}(z)\times
  \Gamma_{q,\,t^{-1}}(qz^{-1})\,
  =\,t^{1/2}\frac{\theta_1\bigl(\,(tz)^{1/2};\,q\bigr)}
  {\theta_1\bigl(z^{1/2};\,q\bigr)}\,,
 \ee
where we take into account  the product representation
 \be\label{EllipticTheta}
  \theta_1(z;\,q)\,
  =\,q^{1/4}\,\frac{z-z^{-1}}{\i}\,\prod_{j\geq1}(1-q^j)(1-z^2q^j)(1-z^{-2}q^j)\,.
 \ee
of the standard elliptic theta-function $\theta_1(z;\,q)$.

Finally, the analog of the Euler integral representation
\eqref{GammaEulerRep} and its inverse \eqref{InversI} are given by
 \be\label{qtGammaEulerRep1}
  \Gamma_{q,t}(x)\,\,
  =\,\sum_{\la\geq 0}\,x^{\la}\,
  \frac{\Gamma_{q,\,tq^{-1}}(q)}{\Gamma_{q,\,tq^{-1}}(q^{\la+1})}\,,
 \hspace{1.2cm}
  \frac{\Gamma_{q,\,tq^{-1}}(q)}{\Gamma_{q,\,tq^{-1}}(q^{\la+1})}\,
  =\,\int_Td^{\times}\!x\,\,x^{-\la}\,
  \Gamma_{q,t}(x)\,.
 \ee

Consider now a specialization of the $\Gamma_{q,t}$ at $t=0$ given
by
 \be\label{qGammaWeierstrass}
  \Gamma_q(z)\,:=\,\frac{1}{(z;\,q)_{\infty}}\,
  =\,\prod_{j=0}^{\infty}\frac{1}{1-zq^j}\,.
 \ee
The $q$-Gamma function $\Gamma_q(z)$ has poles at
$z=q^{-m},\,m\in\IZ_+$ and  satisfy proper analogs of
\eqref{GammaEulerRep}-\eqref{GammaFunctionalEquation}. The
$q$-analog of the Weierstrass product formula
\eqref{GammaWeierstrass} is given by \eqref{qGammaWeierstrass}. The
$q$-analog of the Euler integral formula \eqref{GammaEulerRep} is
given by
 \be
  \Gamma_q(z)\,
  =\,\sum_{\la\geq0}z^{\la}\,\frac{\Gamma_q(q)}{\Gamma_q(q^{\la+1})}\,,
 \hspace{1.5cm}
  \frac{\Gamma_q(q)}{\Gamma_q(q^{\la+1})}\,
  =\,\int_T\!d^{\times}\!z\,\,z^{-\la}\,\Gamma_q(z)\,,
 \ee
and the $q$-analog of the functional equation
\eqref{GammaFunctionalEquation} has the  following form:
 \be
  \Gamma_q(z)\times \Gamma_q(qz^{-1})\,
  =\,\frac{q^{1/4}}{\Gamma_q(q)}
  \frac{\i z^{-1/2}}{\theta_1(z^{1/2};\,q)}\,,
 \ee
which can be deduced  from \eqref{EllipticTheta}.

Now consider the following analog of the $\Gamma$-function
 \be\label{kGammaFunction}
  \Gamma^{(\kappa)}(z)\,
  =\,\lim_{\hbar\to 0}\Gamma_{q,t}(x)\,
  =\,\Big(\frac{1}{1-z}\Big)^{\kappa}\,, \qquad t=e^{\kappa \hbar},
  \quad q=e^{\hbar},
 \ee
depending on a positive integer parameter $\kappa$. The analog of
the functional equation \eqref{GammaFunctionalEquation} reads
 \be
  \Gamma^{(\kappa)}(z)\times\Gamma^{(-\kappa)}(z^{-1})\,
  =\,-\frac{1}{z^{\kappa}}\,,
 \ee
and the binomial formula for
$\Gamma^{(\kappa)}(z)=(1-z)^{-\kappa}$ implies the following analogs
of the Euler's integral formula \eqref{GammaEulerRep} and its inverse:
 \be
  \Gamma^{(\kappa)}(z)\,
  =\,\sum_{n\geq0}\frac{z^n}{n!}\,\frac{\Gamma(\kappa+n)}{\Gamma(\kappa)}\,,
 \qquad
  \frac{\Gamma(\kappa+n)}{\Gamma(\kappa)}\,
  =\,n!\int_T\!d^{\times}\!z\,\,z^n\,\,\Gamma^{(\kappa)}(z)\,.
 \ee
The functions $\Gamma_{q,t}(x)$, $\Gamma_q(x)$ and
$\Gamma^{(\kappa)}(x)$ play an important role in the Baxter
operator formalism for Macdonald, $q$-Whittaker and Jack polynomials
correspondingly.

\vskip 1cm

\noindent {\small {\bf A.G.} {\sl Institute for Theoretical and
Experimental Physics, 117259, Moscow,  Russia; \hspace{8 cm}\,
\hphantom{xxx}  \hspace{2 mm} School of Mathematics, Trinity College
Dublin, Dublin 2, Ireland; \hspace{6 cm}\hspace{5 mm}\,
\hphantom{xxx}   \hspace{2 mm} Hamilton Mathematics Institute,
Trinity College Dublin, Dublin 2, Ireland;}}\\

\noindent{\small {\bf D.L.} {\sl
 Institute for Theoretical and Experimental Physics,
117259, Moscow, Russia};\\
\hphantom{xxxx} {\it E-mail address}: {\tt lebedev@itep.ru}}\\

\noindent{\small {\bf S.O.} {\sl
 Institute for Theoretical and Experimental Physics,
117259, Moscow, Russia};\\
\hphantom{xxxx} {\it E-mail address}: {\tt Sergey.Oblezin@itep.ru}}

\end{document}